\documentclass[12pt,thmsa]{article}
\usepackage{amsfonts}
\usepackage{amsmath}
\usepackage{sw20lart}

\input tcilatex
\begin{document}

\author{Nick Laskin\thanks{\textit{E-mail address}: nlaskin@rocketmail.com}}
\title{\textbf{Fractional non-homogeneous\ counting process}}
\date{TopQuark Inc.\\
Toronto, ON, M6P 2P2\\
Canada}
\maketitle

\begin{abstract}
A new fractional non-homogeneous counting process has been introduced and
developed using the Kilbas and Saigo three-parameter generalization of the
Mittag-Leffler function. The probability distribution function of this
process reproduces for certain set of the fractality parameters the famous
Poisson and fractional Poisson probability distributions as well as the
probability distribution function of a counting process displaying stretched
exponential interarrival times distribution.

Applications of the developed fractional non-homogeneous counting process
cover fractional compound process, the generalization of combinatorial
polynomials and numbers, statistics of the fractional non-homogeneous
counting process and a new representation of the Kilbas and Saigo function.


\end{abstract}

\section{Introduction}

Recently, a discrete probability distribution function

\begin{equation}
P_{\sigma }(n,x)=\frac{x^{\sigma n}}{n!}\exp (-x^{\sigma }),\quad x\geq
0,\quad 0<\sigma \leq 1,\quad n=0,1,2,...,  \label{eq1}
\end{equation}

has been implemented to introduce and study stretched coherent states \cite%
{LaskinArXiV}. Substitution $x\rightarrow \lambda _{\sigma }^{1/\sigma }t$
leads to the probability distribution function of the following counting
process

\begin{equation}
P_{\sigma }(n,t)=\frac{(\lambda _{\sigma }t^{\sigma })^{n}}{n!}\exp
(-\lambda _{\sigma }t^{\sigma }),  \label{eq2}
\end{equation}

where $\lambda _{\sigma }$ has units of $[\lambda _{\sigma }]=time^{-\sigma
}.$ This expression is recognized as a non-homogeneous Poisson process whose
\textquotedblleft rate\textquotedblright\ $r(t)$ changes with time,

\begin{equation}
r(t)=\sigma \lambda _{\sigma }t^{\sigma -1}.  \label{eq.1.10c}
\end{equation}

Then $P_{\sigma }(n,t)$ can be represented as follows,

\begin{equation}
P_{\sigma }(n,t)=\frac{(\Lambda _{\sigma }(t))^{n}}{n!}\exp (-\Lambda
_{\sigma }(t)),  \label{eq3}
\end{equation}

where

\begin{equation}
\Lambda _{\sigma }(t)=\int\limits_{0}^{t}d\tau r(\tau )=\sigma \lambda
_{\sigma }\int\limits_{0}^{t}d\tau \tau ^{\sigma -1}=\lambda _{\sigma
}t^{\sigma }.  \label{eq.1.10e}
\end{equation}

Given the probability distribution $P_{\sigma }(n,t)$, it is easy to obtain
the probability distribution function of interarrival times, i.e., the time
intervals between two consecutive arrivals. It turns out, that for the
counting process described by the probability distribution function
specified by Eq.(\ref{eq2}), interarrival times probability distribution
function $\psi _{\sigma }(\tau )$ has the form

\begin{equation}
\psi _{\sigma }(\tau )=\sigma \lambda _{\sigma }\tau ^{\sigma -1}\exp
(-\lambda _{\sigma }\tau ^{\sigma }),\qquad \tau >0,  \label{eq4}
\end{equation}

which is the well-known Weibull distribution with the shape parameter $%
\sigma $ and the scale parameter $\lambda _{\sigma }$. The probability
distribution function $\psi _{\sigma }(\tau )$ is widely used to fit
\textquotedblleft fat-tail\textquotedblright\ economic, social and natural
data and is called the \textit{stretched exponential distribution }(see, for
example,\textit{\ }\cite{LAHERRERE}).

The function $P_{\sigma }(n,t)$ satisfies the following system of the
equations

\begin{equation}
\frac{\partial P_{\sigma }(n,t)}{\partial t}=\sigma \lambda _{\sigma
}t^{\sigma -1}(P_{\sigma }(n-1,t)-P_{\sigma }(n,t)),\qquad n\geq 1,
\label{eq1.3}
\end{equation}

\begin{equation}
\frac{\partial P_{\sigma }(0,t)}{\partial t}=-\sigma \lambda _{\sigma
}t^{\sigma -1}P_{\sigma }(0,t),  \label{eq1.4}
\end{equation}

with some initial condition at $t=0$,

\begin{equation}
P_{\sigma }(n,t=0)=p_{0}(n).  \label{eq5}
\end{equation}

Usually, it is supposed that at $t=0$ nothing had arrived, that is the
initial condition is specified as $P_{\sigma }(n,t)|_{t=0}=p_{0}(n)$, where $%
p_{0}(n)=0,\ \mathrm{if}\ n\neq 0$ and $p_{0}(n)=1,\ \mathrm{if}\ n=0$.

Thus, Eq.(\ref{eq1.3}) is linear differential equation with respect to $t$
and difference equation with respect to $n$, and is called a
differential-difference equation. When $\sigma =1$, we get the famous
Poisson probability distribution function $P_{1}(n,t)\equiv P_{\sigma
}(n,t)|_{\sigma =1}$, which satisfies the system of equations \cite{Feller}

\begin{equation}
\frac{\partial P_{1}(n,t)}{\partial t}=\lambda
_{1}(P_{1}(n-1,t)-P_{1}(n,t)),\qquad n\geq 1,  \label{eq6a}
\end{equation}

\begin{equation}
\frac{\partial P_{\sigma }(0,t)}{\partial t}=-\lambda _{1}P_{\sigma }(0,t),
\label{eq6b}
\end{equation}

and has the form

\begin{equation}
P_{1}(n,t)=\frac{(\lambda _{1}t)^{n}}{n!}\exp (-\lambda _{1}t),\qquad
n=0,1,2,...,.  \label{eq7}
\end{equation}

where $\lambda _{1}$ has units of $[\lambda _{1}]=time^{-1}.$ As can be seen
from Eq.(\ref{eq4}), when $\sigma =1$, $\psi _{\sigma }(\tau )|_{\sigma =1}=$
$\psi _{1}(\tau )$ becomes an exponential distribution,

\begin{equation}
\psi _{1}(\tau )=\lambda _{1}\exp (-\lambda _{1}\tau ),\qquad \tau >0.
\label{eq8}
\end{equation}

Therefore, $P_{\sigma }(n,t)$ defined by\ Eq.(\ref{eq2}) is generalization
of the Poisson probability distribution function $P_{1}(n,t)$ defined by Eq.(%
\ref{eq7}) and $\psi _{\sigma }(\tau )$ is generalization of the exponential
distribution defined by Eq.(\ref{eq8}).

Another important and widely used generalization of the Poisson process is
fractional Poisson process developed to model and study a long-memory impact
on the counting processes \cite{Laskin7} empirically observed in complex
quantum and classical systems. The quantum system example is the
fluorescence intermittency of single CdSe quantum dots, that is, the
fluorescence emission of single nanocrystals exhibits intermittent behavior,
namely, a sequence of \textquotedblright light on\textquotedblright\ and
\textquotedblright light off\textquotedblright\ states departing from
Poisson statistics. In fact, the waiting time distribution in both states
displays non-exponential behavior \cite{Kuno}. As examples of classical
systems let's mention the distribution of waiting times between two
consecutive transactions in financial markets \cite{Sabatelli} and another,
which comes from network communication systems, where the duration of
network sessions or connections exhibits non-exponential power-law behavior 
\cite{Willinger}. The probability distribution function of the fractional
Poisson process satisfies the following system of fractional
differential-difference equations

\begin{equation}
_{0}^{C}D_{t}^{\mu }P_{\mu }(n,t)=\lambda _{\mu }\left( P_{\mu
}(n-1,t)-P_{\mu }(n,t)\right) ,\qquad n\geq 1,  \label{eqKF1b}
\end{equation}

and

\begin{equation}
_{0}^{C}D_{t}^{\mu }P_{\mu }(0,t)=-\lambda _{\mu }P_{\mu }(0,t),\qquad 0<\mu
\leq 1,  \label{eqKF2b}
\end{equation}

subject to normalization condition

\begin{equation}
\sum\limits_{n=0}^{\infty }P_{\mu }(n,t)=1,  \label{eqKF3}
\end{equation}

and the initial condition at $t=0$, $P_{\mu }(n,t=0)=\delta _{n,0},$ where $%
\mu $ is fractality parameter, the operator $_{0}^{C}D_{t}^{\mu }$ is the
Caputo fractional derivative\footnote{%
The basic formulas on fractional calculus can be found in Refs. \cite{Oldham}
- \cite{Miller}.} whose action on a function $\mathrm{f}(t)$ is defined as 
\cite{Caputo}

\begin{equation}
_{0}^{C}D_{t}^{\mu }\mathrm{f}(t)=\frac{1}{\Gamma (1-\mu )}%
\int\limits_{0}^{t}\frac{d\tau (d\mathrm{f}(\tau )/d\tau )}{(t-\tau )^{\mu }}%
,\qquad 0<\mu \leq 1,  \label{eqKF3b}
\end{equation}

with $\Gamma (\mu )$ being the Gamma function, which has the familiar
representation $\Gamma (\mu )=\int\limits_{0}^{\infty }dte^{-t}t^{\mu -1}$, $%
\mathrm{Re}\mu >0$, and the parameter $\lambda _{\mu }$ has units $[\lambda
_{\mu }]=time^{-\mu }$.

The solution of the system of the equations (\ref{eqKF1b}) and (\ref{eqKF2b}%
) has the form

\begin{equation}
P_{\mu }(n,t)=\frac{(\lambda _{\mu }t^{\mu })^{n}}{n!}\sum\limits_{k=0}^{%
\infty }\frac{(k+n)!}{k!}\frac{(-\lambda _{\mu }t^{\mu })^{k}}{\Gamma (\mu
(k+n)+1))},\quad n=0,1,2,...,\quad 0<\mu \leq 1,  \label{eq9}
\end{equation}

which is probability distribution function of the \textit{fractional Poisson
process} \cite{Laskin7}. The $P_{\mu }(n,t)$ gives us the probability that
in the time interval $[0,t]$ we observe $n$ arrivals.

Equations (\ref{eqKF1b}) and (\ref{eqKF2b}) are fractional generalization of
Eqs.(\ref{eq1.3}) and (\ref{eq1.4}) for probability distribution function of
the Poisson process. When $\mu =1$, Eq.(\ref{eq9}) turns into Eq.(\ref{eq7}).

The goal of this work is to develop a system of fractional non-homogeneous
differential-difference equations the exact solution of which is new
fractional probability distribution of non-homogeneous fractional counting
process. It turns out that this system reproduces for certain set of the
fractality parameters the systems of equations (\ref{eq1.3}) and (\ref{eq1.4}%
), (\ref{eq6a}) and (\ref{eq6b}), (\ref{eqKF1b}) and (\ref{eqKF2b}). This
goal was achieved by implementing the Kilbas and Saigo function \cite%
{Kilbas1}, which is a three-parameter generalization of the Mittag-Leffler
function \cite{ML}.

The paper is organized as follows.

In Section 1 a non-homogeneous fractional differential-difference equation
for the counting probability distribution functions has been developed and
discussed. To solve this equation the probability generating function was
calculated in Section 2. The fractional counting probability distribution
function has been obtained in Section 3 using the probability generating
function. Section 4 is devoted to the calculation of the moment generating
function. The mean and variance of the new fractional counting process were
found using the moment generating function. Interarrival time distribution
was obtained in Section 5. Section 6 presents applications of new fractional
counting probability distribution function covering fractional compound
process and the generalization of combinatorial polynomials and numbers.

The results obtained are summarized in the Conclusion.

The Appendix A lists special cases of Kilbas and Saigo three-parameter
function.

The Appendix B presents analytical expressions for the first four central
moments of the introduced fractional non-homogeneous probability
distribution.

\section{Fundamentals}

We define a new fractional non-homogeneous counting process by introducing
the following system of fractional differential-difference equations for
probability distribution function $P_{\mu ,\beta }(n,t)$ of this process

\begin{equation}
_{0}D_{t}^{\mu }P_{\mu ,\beta }(n,t)=\lambda _{\mu +\beta }t^{\beta }\left(
P_{\mu ,\beta }(n-1,t)-P_{\mu ,\beta }(n,t)\right) ,\quad n=1,2,3,...,\quad
t\geq 0,  \label{eq1.1}
\end{equation}

\begin{equation}
_{0}D_{t}^{\mu }P_{\mu ,\beta }(0,t)=-\lambda _{\mu +\beta }t^{\beta }P_{\mu
,\beta }(0,t),  \label{eq1.1a}
\end{equation}

subject to normalization condition

\begin{equation}
\sum\limits_{n=0}^{\infty }P_{\mu ,\beta }(n,t)=1,  \label{eq1.2}
\end{equation}

and the initial condition

\begin{equation}
P_{\mu ,\beta }(n,t=0)=\delta _{n,0}.  \label{eq1.2a}
\end{equation}

Here, $\mu $ and $\beta $ are dimensionless fractality parameters subject to
the conditions,

\begin{equation}
0<\mu \leq 1,\qquad -\mu <\beta \leq 1-\mu ,  \label{eq1.2ab}
\end{equation}

the rate of arrivals $\lambda _{\mu +\beta }$ has units $[\lambda _{\mu
+\beta }]=$ $time^{-(\mu +\beta )}$, $_{0}D_{t}^{\mu }$ is the operator of
the Caputo fractional time derivative of order $\mu $ defined by Eq.(\ref%
{eqKF3b}). We will see later where the conditions for $\mu $ and $\beta $
come from.

The system of Eqs.(\ref{eq1.1}) and (\ref{eq1.2a}) introduces the counting
process with probability distribution function $P_{\mu ,\beta }(n,t)$ of
arriving $n$ items ($n=0,1,2,...$) by time $t,$ when the rate of arrivals
depends on time $t$. In other words, the above defined counting process is
generalization to non-homogeneous case (due to the dependency of the arrival
rate on time) of the fractional Poisson process originally introduced and
developed in \cite{Laskin7}.

When $\mu =1$ and $\beta =0$, the system of Eqs.(\ref{eq1.1}) and (\ref%
{eq1.2a}) turns into the system of equations (\ref{eq6a}) and (\ref{eq6b}),
which define the Poisson process with rate of arrival $\lambda _{1}$.

When $\beta =0$, the system of Eqs.(\ref{eq1.1}) and (\ref{eq1.2a}) turns
into the system of equations (\ref{eqKF1b}) and (\ref{eqKF2b}), which define
the fractional Poisson process with rate of arrival $\lambda _{\mu }$ where $%
0<\mu \leq 1$.

When $\mu =1$, the system of Eqs.(\ref{eq1.1}) and (\ref{eq1.2a}) turns into
the the system of equations (\ref{eq1.3}) and (\ref{eq1.4}), which define
the probability distribution function $P_{\sigma }(n,t)$, with $\sigma
=1+\beta $, $0<\sigma \leq 1$.

\section{Probability generating function}

To solve the system of equations (\ref{eq1.1}) - (\ref{eq1.2a}) it is
convenient to use the method of the generating function. Introducing the
probability generating function $G_{\mu ,\beta }(s,t)$

\begin{equation}
G_{\mu ,\beta }(s,t)=\sum\limits_{n=0}^{\infty }s^{n}P_{\mu ,\beta }(n,t),
\label{eq1.3ab}
\end{equation}

allows to calculate\ $P_{\mu ,\beta }(n,t)$ as

\begin{equation}
P_{\mu ,\beta }(n,t)=\frac{1}{n!}\frac{\partial ^{n}G_{\mu ,\beta }(s,t)}{%
\partial s^{n}}|_{s=0}.  \label{eq1.4ab}
\end{equation}

Then, by multiplying Eqs.(\ref{eq1.1}) and (\ref{eq1.1a}) by $s^{n}$,
summing over $n$, we obtain the following fractional differential-difference
equation for the probability generating function $G_{\mu }(s,t)$

\begin{equation}
_{0}D_{t}^{\mu }G_{\mu }(s,t)=\lambda _{\mu +\beta }t^{\beta }\left(
\sum\limits_{n=1}^{\infty }s^{n}P_{\mu }(n-1,t)-\sum\limits_{n=0}^{\infty
}s^{n}P_{\mu }(n,t)\right) =  \label{1.5}
\end{equation}

\begin{equation*}
\lambda _{\mu +\beta }t^{\beta }(s-1)G_{\mu }(s,t).
\end{equation*}

Thus, we come to the following equation

\begin{equation}
_{0}D_{t}^{\mu }G_{\mu ,\beta }(s,t)=\lambda _{\mu +\beta }t^{\beta
}(s-1)G_{\mu ,\beta }(s,t),  \label{eq1.8}
\end{equation}

subject to the initial condition

\begin{equation}
G_{\mu ,\beta }(s,t=0)=1,  \label{eq1.9}
\end{equation}

due to the definition Eq.(\ref{eq1.3ab}) and the initial condition for $%
P_{\mu ,\beta }(n,t)$ given by Eq.(\ref{eq1.2a}). The solution to Eq.(\ref%
{eq1.8}) can be obtained using the ansatz \cite{Kilbas1}

\begin{equation}
G_{\mu ,\beta }(s,t)=\sum\limits_{n=0}^{\infty }\frac{c_{n}(\mu ,\beta
)\{\lambda _{\mu +\beta }t^{\mu +\beta }(s-1)\}^{n}}{\Gamma (n(\mu +\beta
)+1)},  \label{eq1.10}
\end{equation}

where $c_{n}(\mu ,\beta )$ are still unknown coefficients and $\Gamma (\mu )$
is the Gamma function. Substituting this ansatz into Eq.(\ref{eq1.8}) and
taking into account that due to the definition of Caputo fractional
derivative Eq.(\ref{eqKF3b})

\begin{equation}
_{0}D_{t}^{\mu }t^{n(\mu +\beta )}=\frac{\Gamma (n(\mu +\beta )+1)}{\Gamma
(n(\mu +\beta )+1-\mu )}t^{n(\mu +\beta )},  \label{eq1.11}
\end{equation}

we find that Eq.(\ref{eq1.8}) is satisfied if the following recurrence
relation holds,

\begin{equation}
c_{n+1}(\mu ,\beta )=\frac{\Gamma (n(\mu +\beta )+\beta +1)}{\Gamma (n(\mu
+\beta )+1)}c_{n}(\mu ,\beta ).  \label{eq1.12}
\end{equation}

This recurrence relation allows us to represent Eq.(\ref{eq1.10}) as

\begin{equation}
G_{\mu ,\beta }(s,t)=c_{0}(\mu ,\beta )+  \label{eq1.13}
\end{equation}

\begin{equation*}
\sum\limits_{n=1}^{\infty }\{\lambda _{\mu +\beta }t^{\mu +\beta
}(s-1)\}^{n}\prod\limits_{k=0}^{n-1}\frac{\Gamma (k(\mu +\beta )+\beta +1)}{%
\Gamma (k(\mu +\beta )+\mu +\beta +1)}c_{0}(\mu ,\beta ).
\end{equation*}

Taking into account the initial condition given by Eq.(\ref{eq1.9}), we come
to the conclusion that $c_{0}(\mu ,\beta )=1$. Therefore, the solution of
the initial problem defined by Eqs.(\ref{eq1.8}) and (\ref{eq1.9}) is

\begin{equation}
G_{\mu ,\beta }(s,t)=1+\sum\limits_{n=1}^{\infty }\{\lambda _{\mu +\beta
}t^{\mu +\beta }(s-1)\}^{n}\prod\limits_{k=0}^{n-1}\frac{\Gamma (k(\mu
+\beta )+\beta +1)}{\Gamma (k(\mu +\beta )+\mu +\beta +1)}.  \label{eq1.14}
\end{equation}

It turns out that the probability generating function $G_{\mu ,\beta }(s,t)$
can be expressed in terms of Kilbas and Saigo three-parameter generalization
of the Mittag-Leffler function \cite{Kilbas1}, defined as

\begin{equation}
E_{\alpha ,m,l}(z)=\sum\limits_{n=0}^{\infty }c_{n}z^{n},\qquad z\in 
\mathbb{C}
,  \label{eq1.15}
\end{equation}

where

\begin{equation}
c_{n}=\prod\limits_{j=0}^{n-1}\frac{\Gamma (\alpha (jm+l)+1)}{\Gamma (\alpha
(jm+l)+1)+1)},\qquad n=1,2,...,\qquad c_{0}=1,  \label{eq1.16}
\end{equation}

and

\begin{equation}
\alpha ,m,l\in 
\mathbb{R}
,\quad \alpha >0,\quad m>0\quad \mathrm{and}\quad \alpha (jm+l)\neq
-1,-2,-3,....  \label{eq1.16a}
\end{equation}

To show how the Kilbas and Saigo function appears in Eq.(\ref{eq1.14}) we
transform the arguments of the Gamma functions under the product sign in Eq.(%
\ref{eq1.14}) as follows

\begin{equation}
\prod\limits_{k=0}^{n-1}\frac{\Gamma (k(\mu +\beta )+\beta +1)}{\Gamma
(k(\mu +\beta )+\mu +\beta +1)}=\prod\limits_{k=0}^{n-1}\frac{\Gamma \lbrack
\mu (k(1+\frac{\beta }{\mu })+\frac{\beta }{\mu })+1]}{\Gamma \lbrack \mu
(k(1+\frac{\beta }{\mu })+\frac{\beta }{\mu }+1)+1]}.  \label{eq1.16b}
\end{equation}

Hence, it can be seen from Eqs.(\ref{eq1.14}) and (\ref{eq1.16b}) that the
power series holds for $G_{\mu ,\beta }(s,t)$

\begin{equation}
G_{\mu ,\beta }(s,t)=1+\sum\limits_{n=1}^{\infty }\{\lambda _{\mu +\beta
}t^{\mu +\beta }(s-1)\}^{n}\prod\limits_{k=0}^{n-1}\frac{\Gamma \lbrack \mu
(k(1+\frac{\beta }{\mu })+\frac{\beta }{\mu })+1]}{\Gamma \lbrack \mu (k(1+%
\frac{\beta }{\mu })+\frac{\beta }{\mu }+1)+1]}.  \label{eq1.16c}
\end{equation}

On the other hand, using the definitions given by Eqs.(\ref{eq1.15}) and (%
\ref{eq1.16}), we see that the right side of Eq.(\ref{eq1.16c}) is the
following Kilbas and Saigo function\footnote{%
The special cases of the function $E_{\mu ,1+\frac{\beta }{\mu },\frac{\beta 
}{\mu }}(z)$ are presented in Appendix A.},

\begin{equation}
E_{\mu ,1+\frac{\beta }{\mu },\frac{\beta }{\mu }}(z)=1+  \label{eq1.16d}
\end{equation}

\begin{equation*}
\sum\limits_{n=1}^{\infty }z^{n}\prod\limits_{k=0}^{n-1}\frac{\Gamma \lbrack
\mu (k(1+\frac{\beta }{\mu })+\frac{\beta }{\mu })+1]}{\Gamma \lbrack \mu
(k(1+\frac{\beta }{\mu })+\frac{\beta }{\mu }+1)+1]},
\end{equation*}

with $z=\lambda _{\mu +\beta }t^{\mu +\beta }(s-1)$. Therefore, the
generating function $G_{\mu ,\beta }(s,t)$ has the form

\begin{equation}
G_{\mu ,\beta }(s,t)=E_{\mu ,1+\frac{\beta }{\mu },\frac{\beta }{\mu }%
}(\lambda _{\mu +\beta }t^{\mu +\beta }(s-1)),  \label{eq1.17}
\end{equation}

where $\mu $ and $\beta $ are subject to the conditions given by Eq.(\ref%
{eq1.2ab}). Thus, we expressed in terms of Kilbas and Saigo three-parameter
function the probability generating function $G_{\mu ,\beta }(s,t)$ of the
counting process introduced by the system of fractional
differential-difference equations (\ref{eq1.1}) and (\ref{eq1.1a}).

When $\mu =1$ and $\beta =0$, then according to Eq.(\ref{eqA1}) the
generating function $G_{1,0}(s,t)$ takes the form

\begin{equation*}
G_{1,0}(s,t)=\exp (\lambda _{1}t(s-1)),
\end{equation*}

which is the generating function of the Poisson probability distribution 
\cite{Feller}.

When $0<\mu \leq 1$ and $\beta =0$, then due to Eq.(\ref{eqA2}) we obtain
from Eq.(\ref{eq1.16c})

\begin{equation}
G_{\mu ,0}(s,t)=1+\sum\limits_{n=1}^{\infty }\frac{\{\lambda _{\mu }t^{\mu
}(s-1)\}^{n}}{\Gamma (\mu n+1)}=E_{\mu }(\lambda _{\mu }t^{\mu }(s-1)),
\label{eq1.17c}
\end{equation}

where $E_{\mu }(z)$ is the Mittag-Leffler function defined as \cite{ML}

\begin{equation}
E_{\mu }(z)=\sum\limits_{n=0}^{\infty }\frac{z^{n}}{\Gamma (\mu n+1)}.
\label{eq1.172d}
\end{equation}

Thus, Eq.(\ref{eq1.17c}) gives us the probability generating function of the
fractional Poisson process \cite{Laskin7}.

When $\mu =1$ and $-1<\beta \leq 0$, then due to Eq.(\ref{eqA3}) and having
introduced the notations,

\begin{equation}
\sigma =1+\beta \qquad \mathrm{and}\qquad \overline{\lambda }_{\sigma }=%
\frac{\lambda _{\sigma }}{\sigma },  \label{eq1.171a}
\end{equation}

we present Eq.(\ref{eq1.16c}) in the form

\begin{equation}
G_{\sigma }(s,t)=\exp (\overline{\lambda }_{\sigma }t^{\sigma }(s-1)),
\label{eq1.172a}
\end{equation}

where $0<\sigma \leq 1$, if $-1<\beta \leq 0$. Thus, Eq.(\ref{eq1.172a}) is
probability generating function of the probability distribution function
given by Eq.(\ref{eq2}).

When $\mu +\beta =1$, then taking into account Eq.(\ref{eqA4}) we can
express Eq.(\ref{eq1.16c}) in the form

\begin{equation*}
G_{\mu ,1-\mu }(s,t)=1+\sum\limits_{n=1}^{\infty }\{\lambda
t(s-1)\}^{n}\prod\limits_{k=0}^{n-1}\frac{\Gamma (k+2-\mu )}{\Gamma (k+2)},
\end{equation*}

or

\begin{equation*}
G_{\mu ,1-\mu }(s,t)=E_{\mu ,\frac{1}{\mu },\frac{1-\mu }{\mu }}(\lambda
t(s-1)),
\end{equation*}

where the following notation was introduced

\begin{equation}
\lambda =\lambda _{\mu +\beta }\left\vert _{\mu +\beta =1}\right. .
\label{eq1.171d}
\end{equation}

\subsection{Probability distribution function}

Using the definition given by Eq.(\ref{eq1.4ab}) we obtain the probability
distribution function $P_{\mu ,\beta }(n,t)$

\begin{equation}
P_{\mu ,\beta }(n,t)=\frac{1}{n!}\frac{\partial ^{n}G_{\mu ,\beta }(s,t)}{%
\partial s^{n}}|_{s=0}=\frac{(\lambda _{\mu +\beta }t^{\mu +\beta })^{n}}{n!}%
E_{\mu ,1+\frac{\beta }{\mu },\frac{\beta }{\mu }}^{(n)}(-\lambda _{\mu
+\beta }t^{\mu +\beta }),  \label{eq1.18}
\end{equation}

where the following notation is introduced

\begin{equation}
E_{\mu ,1+\frac{\beta }{\mu },\frac{\beta }{\mu }}^{(n)}(z)=\frac{d^{n}}{%
dz^{n}}E_{\mu ,1+\frac{\beta }{\mu },\frac{\beta }{\mu }}(z).  \label{eq1.19}
\end{equation}

The equation (\ref{eq1.18}) can be rewritten as

\begin{equation}
P_{\mu ,\beta }(n,t)=\frac{(-z)^{n}}{n!}E_{\mu ,1+\frac{\beta }{\mu },\frac{%
\beta }{\mu }}^{(n)}(z)\left\vert _{z=-\lambda _{\mu +\beta }t^{\mu +\beta
}}\right. ,  \label{eq1.20}
\end{equation}

and

\begin{equation}
P_{\mu ,\beta }(n=0,t)=E_{\mu ,1+\frac{\beta }{\mu },\frac{\beta }{\mu }%
}(-\lambda _{\mu +\beta }t^{\mu +\beta }).  \label{eq1.21}
\end{equation}

The function $P_{\mu ,\beta }(n,t)$ is normalized

\begin{equation}
\sum\limits_{n=0}^{\infty }P_{\mu ,\beta }(n,t)=\frac{(-z)^{n}}{n!}E_{\mu ,1+%
\frac{\beta }{\mu },\frac{\beta }{\mu }}^{(n)}(z)|_{z=-\lambda _{\mu +\beta
}t^{\mu +\beta }}=E_{\mu ,1+\frac{\beta }{\mu },\frac{\beta }{\mu }%
}^{(n)}(0)=1,  \label{eq1.21a}
\end{equation}

due to the definition given by Eqs.(\ref{eq1.15}) and (\ref{eq1.16}), and
positive due to the complete monotonicity of the three-parameter Kilbas and
Saigo function $E_{\mu ,1+\frac{\beta }{\mu },\frac{\beta }{\mu }}(-x)$ when 
$\mu $ and $\beta $ are subject to the conditions\footnote{%
The conditions 
\begin{equation}
0<\mu \leq 1,\qquad -\mu <\beta ,  \label{eqCM1}
\end{equation}%
\par
follow directly from the condition of complete monotonicity of the
three-parameter Kilbas and Saigo function \cite{Boudabsa}. \ We have added
an additional condition
\par
\begin{equation*}
\beta \leq 1-\mu ,
\end{equation*}%
\par
to ensure consistency with the parameter $\sigma $ defined by Eq.(\ref%
{eq1.171a}) and included in Eq.(\ref{eq1}).} defined by Eq.(\ref{eq1.2ab}),
see, \cite{Boudabsa}, \cite{BoudabsaSimon},

\begin{equation}
(-1)^{n}\frac{d^{n}E_{\mu ,1+\frac{\beta }{\mu },\frac{\beta }{\mu }}(-x)}{%
dx^{n}}\geq 0,\qquad x\geq 0,  \label{eq1.22}
\end{equation}

Therefore, $P_{\mu ,\beta }(n,t)$ is indeed a probability distribution
function. The power series representation for $P_{\mu ,\beta }(n,t)$ is

\begin{equation}
P_{\mu ,\beta }(n,t)=\frac{(\lambda _{\mu +\beta }t^{\mu +\beta })^{n}}{n!}%
\times  \label{eq1.24}
\end{equation}

\begin{equation*}
\sum\limits_{m=0}^{\infty }\frac{(m+n)!}{m!}(-\lambda _{\mu +\beta }t^{\mu
+\beta })^{m}\prod\limits_{k=0}^{m+n-1}\frac{\Gamma \lbrack \mu (k(1+\frac{%
\beta }{\mu })+\frac{\beta }{\mu })+1]}{\Gamma \lbrack \mu (k(1+\frac{\beta 
}{\mu })+\frac{\beta }{\mu }+1)+1]}.
\end{equation*}

Replacing $\lambda _{\mu +\beta }t^{\mu +\beta }$ with $x^{\mu +\beta }$ in
Eqs.(\ref{eq1.20}) and (\ref{eq1.24}), we obtain a new family of fractional
probability distributions $P_{\mu ,\beta }(x)$,

\begin{equation}
P_{\mu ,\beta }(n,x)=\frac{(-z)^{n}}{n!}E_{\mu ,1+\frac{\beta }{\mu },\frac{%
\beta }{\mu }}^{(n)}(z)|_{z=-\lambda _{\mu +\beta }t^{\mu +\beta }},
\label{eq1.25}
\end{equation}

or in a series form

\begin{equation}
P_{\mu ,\beta }(n,x)=\frac{x^{(\mu +\beta )n}}{n!}\sum\limits_{m=0}^{\infty }%
\frac{(m+n)!}{m!}\times  \label{eq1.26}
\end{equation}

\begin{equation*}
(-x^{(\mu +\beta )})^{m}\prod\limits_{k=0}^{m+n-1}\frac{\Gamma \lbrack \mu
(k(1+\frac{\beta }{\mu })+\frac{\beta }{\mu })+1]}{\Gamma \lbrack \mu (k(1+%
\frac{\beta }{\mu })+\frac{\beta }{\mu }+1)+1]},
\end{equation*}

where $\mu $ and $\beta $ are subject to the conditions given by Eq.(\ref%
{eq1.2ab}).

When $\mu =1$ and $\beta =0$, using Eq.(\ref{eqA1}) we get from Eq.(\ref%
{eq1.18})

\begin{equation*}
P_{1,0}(n,t)=\frac{(\lambda _{1}t)^{n}}{n!}\exp (-\lambda _{1}t),
\end{equation*}

which is the celebrated Poisson probability distribution.

When $0<\mu \leq 1$ and $\beta =0$, using Eq.(\ref{eqA2}) we get from Eq.(%
\ref{eq1.18})

\begin{equation*}
P_{\mu ,0}(n,t)=\frac{(\lambda _{\mu }t^{\mu })^{n}}{n!}\sum\limits_{m=0}^{%
\infty }\frac{(m+n)!}{m!}\frac{(-\lambda _{\mu }t^{\mu })^{m}}{\Gamma (\mu
m+1)},
\end{equation*}

which is the fractional Poisson probability distribution \cite{Laskin7}.

When $\mu =1$ and $-1<\beta \leq 0$, then according to Eq.(\ref{eqA3}) we
obtain from Eq.(\ref{eq1.18})

\begin{equation}
P_{1,\beta }(n,t)\equiv P_{1+\beta }(n,t)=\frac{(\frac{\lambda _{1+\beta }}{%
1+\beta }t^{1+\beta })^{n}}{n!}\exp (-\frac{\lambda _{1+\beta }t^{1+\beta }}{%
1+\beta }).  \label{eq1.21e1}
\end{equation}

Using the notations given by Eq.(\ref{eq1.171a}), we represent Eq.(\ref%
{eq1.21e1}) in the form

\begin{equation*}
P_{\sigma }(n,t)=\frac{(\overline{\lambda }_{\sigma }t^{\sigma })^{n}}{n!}%
\exp (-\overline{\lambda }_{\sigma }t^{\sigma }),\quad 0<\sigma \leq 1,
\end{equation*}

which is recognized as Eq.(\ref{eq1}).

With $\mu +\beta =1$ and using Eq.(\ref{eqA4}) we can present Eq.(\ref%
{eq1.18}) as

\begin{equation*}
P_{\mu ,1-\mu }(n,t)=\frac{(\lambda t)^{n}}{n!}E_{\mu ,\frac{1}{\mu },\frac{%
1-\mu }{\mu }}^{(n)}(-\lambda t),
\end{equation*}

where $\lambda $ is defined by Eq.(\ref{eq1.171d}). The series
representation for $P_{\mu ,1-\mu }(n,t)$ has the form

\begin{equation}
P_{\mu ,1-\mu }(n,t)=\frac{(\lambda t)^{n}}{n!}\sum\limits_{m=0}^{\infty }%
\frac{(m+n)!}{m!}(-\lambda t)^{m}\prod\limits_{k=0}^{m+n-1}\frac{\Gamma
(k+2-\mu )}{\Gamma (k+2)}.  \label{eq1.122a}
\end{equation}

\section{Moment generating function}

An equation for the moment of any integer order $m$ of the probability
distribution function $P_{\mu ,\beta }(n,t)$ introduced by Eq. (\ref{eq1.20}%
) can be easily found using the moment generating function $M_{\mu ,\beta
}(s,t)$ defined as

\begin{equation}
M_{\mu ,\beta }(s,t)=\sum\limits_{n=0}^{\infty }e^{-sn}P_{\mu ,\beta }(n,t).
\label{eq1.27}
\end{equation}

Then for the moment of $m^{\mathrm{th}}$ order we have

\begin{equation}
<n_{\mu ,\beta }^{m}>=\sum\limits_{n=0}^{\infty }n^{m}P_{\mu ,\beta
}(n,t)=(-1)^{m}\frac{\partial ^{m}}{\partial s^{m}}M_{\mu ,\beta
}(s,t)\left\vert _{s=0}\right. .  \label{eq1.28}
\end{equation}

Multiplying Eq. (\ref{eq1.20}) by $e^{-sn}$ and summing over the $n$ we get

\begin{equation}
M_{\mu ,\beta }(s,t)=E_{\mu ,1+\frac{\beta }{\mu },\frac{\beta }{\mu }%
}(\lambda _{\mu +\beta }t^{\mu +\beta }(e^{-s}-1)),  \label{eq1.29}
\end{equation}

or in a series form

\begin{equation}
M_{\mu ,\beta }(s,t)=1+\sum\limits_{n=1}^{\infty }\{\lambda _{\mu +\beta
}t^{\mu +\beta }(e^{-s}-1)\}^{n}\prod\limits_{k=0}^{n-1}\frac{\Gamma \lbrack
\mu (k(1+\frac{\beta }{\mu })+\frac{\beta }{\mu })+1]}{\Gamma \lbrack \mu
(k(1+\frac{\beta }{\mu })+\frac{\beta }{\mu }+1)+1]}.  \label{eq1.30}
\end{equation}

Replacing $\lambda _{\mu +\beta }t^{\mu +\beta }$ with $x^{\mu +\beta }$ in
Eqs.(\ref{eq1.29}) and (\ref{eq1.30}) we get the moment generating function
for the fractional generalized probability distributions $P_{\mu ,\beta }(x)$%
,

\begin{equation}
M_{\mu ,\beta }(s,x)=E_{\mu ,1+\frac{\beta }{\mu },\frac{\beta }{\mu }%
}(x^{\mu +\beta }(e^{-s}-1)),  \label{eq1.31}
\end{equation}

and in series form

\begin{equation}
M_{\mu ,\beta }(s,x)=1+\sum\limits_{n=1}^{\infty }\{x^{\mu +\beta
}(e^{-s}-1)\}^{n}\prod\limits_{k=0}^{n-1}\frac{\Gamma \lbrack \mu (k(1+\frac{%
\beta }{\mu })+\frac{\beta }{\mu })+1]}{\Gamma \lbrack \mu (k(1+\frac{\beta 
}{\mu })+\frac{\beta }{\mu }+1)+1]}.  \label{eq1.32}
\end{equation}

When $\mu =1$ and $\beta =0$, using Eq.(\ref{eqA1}) we get from Eq.(\ref%
{eq1.32})

\begin{equation*}
M_{1,0}(s,x)=1+\sum\limits_{n=1}^{\infty }\frac{(\lambda _{1}t(e^{-s}-1))^{n}%
}{n!}=\exp (\lambda _{1}t(e^{-s}-1)),
\end{equation*}

which is the moment generating function of the Poisson process.

With $0<\mu \leq 1$ and $\beta =0$, using Eq.(\ref{eqA2}) we get from Eq.(%
\ref{eq1.32})

\begin{equation*}
M_{\mu ,0}(s,x)=E_{\mu }(\lambda _{\mu }t^{\mu }(e^{-s}-1)),
\end{equation*}

which is the moment generating function of the fractional Poisson process 
\cite{Laskin7}.

For $\mu =1$ and $-1<\beta \leq 0$, Eq.(\ref{eq1.32}) gives us the moment
generating function of the counting process, the probability distribution
function of which is determined by Eq.(\ref{eq2})

\begin{equation}
M_{\sigma }(s,x)\equiv M_{1+\beta }(s,t)=\exp (\overline{\lambda }_{\sigma
}t^{\sigma }(e^{-s}-1)),  \label{eq1.32a}
\end{equation}

where $\sigma $ and $\overline{\lambda }_{\sigma }$ are defined by Eq.(\ref%
{eq1.171a}).

When $\mu +\beta =1$, Eq.(\ref{eq1.31}) goes to

\begin{equation*}
M_{\mu ,1-\mu }(s,t)=E_{\mu ,\frac{1}{\mu },\frac{1-\mu }{\mu }}(\lambda
t(e^{-s}-1),
\end{equation*}

and Eq.(\ref{eq1.32}) goes to

\begin{equation*}
M_{\mu ,1-\mu }(s,t)=1+\sum\limits_{n=1}^{\infty }\{\lambda
t(e^{-s}-1)\}^{n}\prod\limits_{k=0}^{n-1}\frac{\Gamma (k+2-\mu )}{\Gamma
(k+2)},
\end{equation*}

where $\lambda $ is defined by Eq.(\ref{eq1.171d}).

\subsection{Mean and variance}

The mean $<n_{\mu ,\beta }>$ can be calculated using the definition (\ref%
{eq1.28}),

\begin{equation}
<n_{\mu ,\beta }>=-\frac{\partial }{\partial s}M_{\mu ,\beta
}(s,t)\left\vert _{s=0}\right. ,  \label{eq1.33}
\end{equation}

with $M_{\mu ,\beta }(s,t)$ given either by Eq.(\ref{eq1.29}) or Eq.(\ref%
{eq1.30}). Thus, we have

\begin{equation}
<n_{\mu ,\beta }>=\frac{\Gamma (\beta +1)}{\Gamma (\mu +\beta +1)}\lambda
_{\mu +\beta }t^{\mu +\beta }.  \label{eq1.34}
\end{equation}

It is easy to see that the moment of the second order $<n_{\mu ,\beta }^{2}>$
is

\begin{equation}
<n_{\mu ,\beta }^{2}>=\frac{\partial ^{2}}{\partial s^{2}}M_{\mu ,\beta
}(s,t)\left\vert _{s=0}\right. =  \label{eq1.35}
\end{equation}

\begin{equation*}
2\frac{\Gamma (\beta +1)\Gamma (\mu +2\beta +1)}{\Gamma (\mu +\beta
+1)\Gamma (2\mu +2\beta +1)}(\lambda _{\mu +\beta }t^{\mu +\beta })^{2}+%
\frac{\Gamma (\beta +1)}{\Gamma (\mu +\beta +1)}\lambda _{\mu +\beta }t^{\mu
+\beta }.
\end{equation*}

Then variance $Var_{\mu ,\beta }$ of the probability distribution $P_{\mu
,\beta }(n,t)$ is

\begin{equation}
Var_{\mu ,\beta }=<n_{\mu ,\beta }^{2}>-<n_{\mu ,\beta }>^{2}=
\label{eq1.36}
\end{equation}

\begin{equation*}
\left\{ \frac{2\Gamma (\beta +1)\Gamma (\mu +2\beta +1)}{\Gamma (\mu +\beta
+1)\Gamma (2\mu +2\beta +1)}-\left( \frac{\Gamma (\beta +1)}{\Gamma (\mu
+\beta +1)}\right) ^{2}\right\} (\lambda _{\mu +\beta }t^{\mu +\beta })^{2}+
\end{equation*}

\begin{equation*}
\frac{\Gamma (\beta +1)}{\Gamma (\mu +\beta +1)}\lambda _{\mu +\beta }t^{\mu
+\beta },
\end{equation*}

here $\mu $ and $\beta $ are subject to the conditions given by Eq.(\ref%
{eq1.2ab}).

\section{Interarrival times distribution}

An important characteristic of any counting process is the time between
successive counts/arrivals - the interarrival time. Interarrival time is a
random variable having the probability distribution function $\psi _{\mu
,\beta }(\tau ),$

\begin{equation}
\psi _{\mu ,\beta }(\tau )=-\frac{d}{d\tau }\mathcal{P}_{\mu ,\beta }(\tau ),
\label{eq1.37}
\end{equation}

where $\mathcal{P}_{\mu ,\beta }(\tau )$ is the probability that a given
interarrival time is greater or equal to $\tau $ \cite{Laskin7}

\begin{equation}
\mathcal{P}_{\mu ,\beta }(\tau )=1-\sum\limits_{n=1}^{\infty }P_{\mu ,\beta
}(n,\tau )=P_{\mu ,\beta }(0,\tau ),  \label{eq1.38}
\end{equation}

here $P_{\mu ,\beta }(n,\tau )$ is defined by Eq.(\ref{eq1.24}).

It follows from Eqs.(\ref{eq1.37}) and (\ref{eq1.38}) that the probability
distribution function of the interarrival times has the form

\begin{equation}
\psi _{\mu ,\beta }(\tau )=\frac{\mu +\beta }{\tau }\sum\limits_{l=0}^{%
\infty }(-1)^{l}(l+1)(\lambda _{\mu +\beta }\tau ^{\mu +\beta
})^{l+1}\prod\limits_{k=0}^{l}\frac{\Gamma \lbrack \mu (k(1+\frac{\beta }{%
\mu })+\frac{\beta }{\mu })+1]}{\Gamma \lbrack \mu (k(1+\frac{\beta }{\mu })+%
\frac{\beta }{\mu }+1)+1]}.  \label{eq1.39}
\end{equation}

The Laplace transform $\varphi _{\mu ,\beta }(u)$ of the interarrival
probability distribution function is defined as

\begin{equation}
\varphi _{\mu ,\beta }(u)=\int\limits_{0}^{\infty }d\tau e^{-u\tau }\psi
_{\mu ,\beta }(\tau ),  \label{eq1.40}
\end{equation}

and has the form

\begin{equation}
\varphi _{\mu ,\beta }(u)=(\mu +\beta )\sum\limits_{l=0}^{\infty
}(-1)^{l}(l+1)(\lambda _{\mu +\beta })^{l+1}\frac{\Gamma \lbrack (\mu +\beta
)(l+1)]}{u^{(\mu +\beta )(l+1)}}\times  \label{eq1.39b}
\end{equation}

\begin{equation*}
\prod\limits_{k=0}^{l}\frac{\Gamma \lbrack \mu (k(1+\frac{\beta }{\mu })+%
\frac{\beta }{\mu })+1]}{\Gamma \lbrack \mu (k(1+\frac{\beta }{\mu })+\frac{%
\beta }{\mu }+1)+1]}.
\end{equation*}

It easy to see that when $\mu =1$ and $\beta =0$, Eq.(\ref{eq1.39})
transforms into the exponential distribution defined by Eq.(\ref{eq8}). From
Eq.(\ref{eq1.39b}) we obtain the Laplace transform of the exponential
distribution

\begin{equation*}
\varphi _{1,0}(u)=\frac{\lambda _{1}}{\lambda _{1}+u}.
\end{equation*}

When $0<\mu \leq 1$ and $\beta =0$, Eq.(\ref{eq1.39}) becomes the
probability distribution of interarrival times of the fractional Poisson
process

\begin{equation}
\psi _{\mu ,0}(\tau )=\lambda _{\mu }\tau ^{\mu -1}E_{\mu ,\mu }(-\lambda
_{\mu }\tau ^{\mu }),  \label{eq1.39c}
\end{equation}

where $E_{\mu ,\mu }(z)$ is the two parameter Mittag-Leffler function
defined by the series \cite{Erdelyi}

\begin{equation}
E_{\mu ,\mu }(z)=\sum\limits_{n=0}^{\infty }\frac{z^{n}}{\Gamma (\mu n+\mu )}%
.  \label{eq1.39d}
\end{equation}

From and Eq.(\ref{eq1.39b}) it follows that the Laplace transform has the
form \cite{Laskin7}

\begin{equation*}
\varphi _{\mu ,0}(u)=\frac{\lambda _{\mu }}{\lambda _{\mu }+u^{\mu }}.
\end{equation*}

When $\mu =1$ and $-1<\beta \leq 0$, Eq.(\ref{eq1.39}) goes to stretched
exponential distribution defined by Eq.(\ref{eq4}) if we take into account
the notation introduced by Eq.(\ref{eq1.171a}). From Eq.(\ref{eq1.39b}) we
obtain the Laplace transform of the stretched exponential distribution

\begin{equation}
\varphi _{\sigma }(u)=\sigma \overline{\lambda }_{\sigma
}\sum\limits_{l=0}^{\infty }\frac{(-\overline{\lambda }_{\sigma })^{l}}{l!}%
\frac{\Gamma (\sigma (l+1))}{u^{\sigma (l+1)}},  \label{eq1.40a}
\end{equation}

here $\sigma $ and $\overline{\lambda }_{\sigma }$ are defined by Eq.(\ref%
{eq1.171a}).

With $\mu +\beta =1$ Eq.(\ref{eq1.39}) goes to

\begin{equation}
\psi _{\mu ,1-\mu }(\tau )=\frac{1}{\tau }\sum\limits_{l=0}^{\infty
}(-1)^{l}(l+1)(\lambda \tau )^{l+1}\prod\limits_{k=0}^{l}\frac{\Gamma
(k+2-\mu )}{\Gamma (k+2)},  \label{eq1.40b}
\end{equation}

and Eq.(\ref{eq1.39b}) goes to

\begin{equation}
\varphi _{\mu ,1-\mu }(u)=\frac{\lambda }{u}\sum\limits_{l=0}^{\infty }(-%
\frac{\lambda }{u})^{l}(l+1)!\prod\limits_{k=0}^{l}\frac{\Gamma (k+2-\mu )}{%
\Gamma (k+2)},  \label{eq1.122ab}
\end{equation}

where $\lambda $ is defined by Eq.(\ref{eq1.171d}).

\section{Applications}

\subsection{Fractional compound process}

We call stochastic process $\{X_{\mu ,\beta }(t)$, $t\geq 0\}$ a fractional
compound process if it is represented by

\begin{equation}
X_{\mu ,\beta }(t)=\sum\limits_{i=1}^{N_{\mu ,\beta }(t)}Y_{i},
\label{eq4.1}
\end{equation}

where $\{N_{\mu ,\beta }(t)$, $t\geq 0\}$ is a fractional counting process
whose probability distribution function is introduced by Eq.(\ref{eq1.20})
or Eq.(\ref{eq1.24}), and $\{Y_{i}$, $i=1,2,...\}$ are independent and
identically distributed random variables with probability distribution
function $p(Y_{i})$ for each $Y_{i}$. It is assumed that the process $%
\{N_{\mu ,\beta }(t)$, $t\geq 0\}$ and the sequence of random variables $%
\{Y_{i}$, $i=1,2,...\}$ are independent.

The moment generating function $J_{\mu ,\beta }(s,t)$ of fractional
generalized compound process is defined as follows

\begin{equation}
J_{\mu ,\beta }(s,t)=<\exp \{sX_{\mu ,\beta }(t)\}>_{Y,N_{\mu ,\beta }(t)},
\label{eq4.2}
\end{equation}

where $<...>_{Y_{i},N_{\mu ,\beta }(t)}$ stands for average which includes
two statistically independent averaging procedures:

1. Averaging $<...>_{Y}$ over independent random variables $Y_{i}$, defined
by

\begin{equation}
<...>_{Y}=\int
dY_{1}...dY_{n}p(Y_{1})...p(Y_{n})...=\prod\limits_{i=1}^{n}\int
dY_{i}p(Y_{i})....  \label{eq4.3}
\end{equation}

2. Averaging over random number of counts defined by

\begin{equation}
<...>_{N_{\mu ,\beta }(t)}=\sum\limits_{n=0}^{\infty }P_{\mu ,\beta
}(n,t)...,  \label{eq4.4}
\end{equation}

where $P_{\mu ,\beta }(n,t)$ is given by Eq.(\ref{eq1.20}).

One can see from Eq.(\ref{eq4.2}) that $k^{\mathrm{th}}$ order moment of
fractional compound process $X_{\mu ,\beta }(t)$ is obtained by
differentiating $J_{\mu ,\beta }(s,t)$ $k$ times with respect to $s$ and
then setting $s=0$, that is

\begin{equation}
<(X_{\mu ,\beta }(t))^{k}>_{Y,N_{\mu ,\beta }(t)}=\frac{\partial ^{k}}{%
\partial s^{k}}J_{\mu ,\beta }(s,t)_{|s=0}.  \label{eq4.5}
\end{equation}

To obtain equation for the moment generating function $J_{\mu ,\beta }(s,t)$
we apply Eqs.(\ref{eq4.3}) and (\ref{eq4.4}) to Eq.(\ref{eq4.2})

\begin{equation}
J_{\mu ,\beta }(s,t)=\sum\limits_{n=0}^{\infty }<\exp \{sX_{\mu ,\beta
}(t)|N_{\mu ,\beta }(t)=n\}>_{Y}P_{\mu ,\beta }(n,t)=  \label{eq4.6}
\end{equation}

\begin{equation*}
\sum\limits_{n=0}^{\infty }<\exp \{s(Y_{1}+...+Y_{n})\}>_{Y}\times \frac{%
(-z)^{n}}{n!}E_{\mu ,1+\frac{\beta }{\mu },\frac{\beta }{\mu }%
}^{(n)}(z)\left\vert _{z=-\lambda _{\mu +\beta }t^{\mu +\beta }}\right. ,
\end{equation*}

where we used the independence of $\{Y_{1},Y_{2},...\}$ and $\{N_{\mu ,\beta
}(t)$, $t\geq 0\}$ and the independence of the $Y_{i}$'s between themselves.
Hence, setting

\begin{equation}
g(s)=\int dY_{i}p(Y_{i})\exp \{sY_{i}\},  \label{eq4.7}
\end{equation}

for the moment generating function of random variables $Y_{i}$, we obtain
from Eq.(\ref{eq4.6}) the moment generating function $J_{\mu ,\beta }(s,t)$
of the fractional compound process

\begin{equation}
J_{\mu ,\beta }(s,t)=\sum\limits_{n=0}^{\infty }g^{n}(s)\times \frac{(-z)^{n}%
}{n!}\frac{d^{n}}{dz^{n}}\left( E_{\mu ,1+\frac{\beta }{\mu },\frac{\beta }{%
\mu }}(z)\right) |_{z=-\lambda _{\mu +\beta }t^{\mu +\beta }}=
\label{eq.4.8}
\end{equation}%
\begin{equation*}
E_{\mu ,1+\frac{\beta }{\mu },\frac{\beta }{\mu }}(\lambda _{\mu +\beta
}t^{\mu +\beta }(g(s)-1)).
\end{equation*}

By differentiating the above with respect to $s$ and setting $s=0$, it is
easy to obtain the mean of the fractional compound process,

\begin{equation}
<X_{\mu ,\beta }(t)>_{Y_{i},N_{\mu ,\beta }(t)}=\frac{\partial }{\partial s}%
J_{\mu ,\beta }(s,t)_{|s=0}=<Y>_{Y}\lambda _{\mu +\beta }t^{\mu +\beta }%
\frac{\Gamma (\beta +1)}{\Gamma (\mu +\beta +1)},  \label{eq4.9}
\end{equation}

which is a manifestation of independency between fractional counting process 
$\{N_{\mu ,\beta }(t)$, $t\geq 0\}$ and random variables $Y_{i}.$

\subsection{New polynomials and associated numbers}

Based on the probability distribution function defined by Eq.(\ref{eq1.25})
we introduce polynomials $B_{\mu ,\beta }(x,m)$,

\begin{equation}
B_{\mu ,\beta }(x,m)=\sum\limits_{n=0}^{\infty }n^{m}P_{\mu ,\beta
}(n,x)=\sum\limits_{n=0}^{\infty }n^{m}\frac{x^{n}}{n!}E_{\mu ,1+\frac{\beta 
}{\mu },\frac{\beta }{\mu }}^{(n)}(-x),  \label{eq5.1a}
\end{equation}

\begin{equation*}
B_{\mu ,\beta }(x,0)=1,
\end{equation*}

where the parameters $\mu $, and $\beta $ satisfy the conditions given by
Eq.(\ref{eq1.2ab}). From Eq.(\ref{eq5.1a}) at $x=1$ we come to fractional
numbers $B_{\mu ,\beta }(m)$,

\begin{equation}
B_{\mu ,\beta }(m)=B_{\mu ,\beta }(x,m)|_{x=1}=\sum\limits_{n=0}^{\infty }%
\frac{n^{m}}{n!}E_{\mu ,1+\frac{\beta }{\mu },\frac{\beta }{\mu }}^{(n)}(-1).
\label{eq5.5}
\end{equation}

As an example, here are the first two newly introduced polynomials

\begin{equation}
B_{\mu ,\beta }(x,1)=\frac{\Gamma (\beta +1)}{\Gamma (\mu +\beta +1)}x,
\label{eq5.5a}
\end{equation}

and

\begin{equation}
B_{\mu ,\beta }(x,2)=\frac{\Gamma (\beta +1)\Gamma (\mu +2\beta +1)}{\Gamma
(\mu +\beta +1)\Gamma (2\mu +2\beta +1)}x^{2}+\frac{\Gamma (\beta +1)}{%
\Gamma (\mu +\beta +1)}x.  \label{eq5.5b}
\end{equation}

The first two fractional numbers associated with the polynomials $B_{\mu
,\beta }(x,m)$ are

\begin{equation}
B_{\mu ,\beta }(1)=\frac{\Gamma (\beta +1)}{\Gamma (\mu +\beta +1)},
\label{eq5.5c}
\end{equation}

and

\begin{equation}
B_{\mu ,\beta }(2)=\frac{\Gamma (\beta +1)\Gamma (\mu +2\beta +1)}{\Gamma
(\mu +\beta +1)\Gamma (2\mu +2\beta +1)}+\frac{\Gamma (\beta +1)}{\Gamma
(\mu +\beta +1)}.  \label{eq5.5d}
\end{equation}

Series representations of the polynomials $B_{\mu ,\beta }(x,m)$ and numbers 
$B_{\mu ,\beta }(m)$ are

\begin{equation}
B_{\mu ,\beta }(x,m)=\sum\limits_{n=0}^{\infty }n^{m}\frac{x^{n}}{n!}%
\sum\limits_{r=0}^{\infty }\frac{(r+n)!}{r!}(-x)^{r}\prod%
\limits_{k=0}^{r+n-1}\frac{\Gamma \lbrack \mu (k(1+\frac{\beta }{\mu })+%
\frac{\beta }{\mu })+1]}{\Gamma \lbrack \mu (k(1+\frac{\beta }{\mu })+\frac{%
\beta }{\mu }+1)+1]},  \label{eq5.7}
\end{equation}

\begin{equation}
B_{\mu ,\beta }(m)=B_{\mu ,\beta }(x,m)|_{x=1}=\sum\limits_{n=0}^{\infty }%
\frac{n^{m}}{n!}\times  \label{eq5.8}
\end{equation}

\begin{equation*}
\sum\limits_{r=0}^{\infty }\frac{(r+n)!}{r!}\prod\limits_{k=0}^{r+n-1}\frac{%
\Gamma \lbrack \mu (k(1+\frac{\beta }{\mu })+\frac{\beta }{\mu })+1]}{\Gamma
\lbrack \mu (k(1+\frac{\beta }{\mu })+\frac{\beta }{\mu }+1)+1]},
\end{equation*}

where we used Eq.(\ref{eq1.26}).

\subsubsection{Generating function of polynomials $B_{\protect\mu ,\protect%
\beta }(x,m)$}

The generating function of the polynomials $B_{\mu ,\beta }(x,m)$ is defined
as

\begin{equation}
\mathcal{B}_{\mu ,\beta }(s,x)=\sum\limits_{m=0}^{\infty }\frac{s^{m}}{m!}%
B_{\mu ,\beta }(x,m).  \label{eq6.1}
\end{equation}

To get the polynomial $B_{\mu ,\beta }(x,m)$ we differentiate $\mathcal{B}%
_{\mu ,\beta }(s,x)$ $m$ times with respect to $s$, and then let $s=0$. That
is,

\begin{equation}
B_{\mu ,\beta }(x,m)=\frac{\partial ^{m}}{\partial s^{m}}\mathcal{B}_{\mu
,\beta }(s,x)|_{s=0}.  \label{eq6.2}
\end{equation}

To find an explicit equation for $\mathcal{B}_{\mu ,\beta }(s,x)$, let's
substitute Eq.(\ref{eq5.1a}) into Eq.(\ref{eq6.1})

\begin{equation*}
\mathcal{B}_{\mu ,\beta }(s,x)=\sum\limits_{m=0}^{\infty }\frac{s^{m}}{m!}%
\sum\limits_{n=0}^{\infty }n^{m}\frac{x^{n}}{n!}E_{\mu ,1+\frac{\beta }{\mu }%
,\frac{\beta }{\mu }}^{(n)}(-x).
\end{equation*}

Summation over $m$ gives us

\begin{equation*}
\sum\limits_{m=0}^{\infty }\frac{s^{m}}{m!}n^{m}=e^{sn},
\end{equation*}

and we obtain

\begin{equation}
\mathcal{B}_{\mu ,\beta }(s,x)=\sum\limits_{n=0}^{\infty }\frac{x^{n}}{n!}%
e^{sn}E_{\mu ,1+\frac{\beta }{\mu },\frac{\beta }{\mu }}^{(n)}(-x)=E_{\mu ,1+%
\frac{\beta }{\mu },\frac{\beta }{\mu }}(x(e^{s}-1)),  \label{eq6.3}
\end{equation}

where the function $E_{\mu ,1+\frac{\beta }{\mu },\frac{\beta }{\mu }}(z)$
is defined by Eq.(\ref{eq1.16d}). It is easy to see that the representation
of the generating function $\mathcal{B}_{\mu ,\beta }(s,x)$ in a power
series has the form

\begin{equation}
\mathcal{B}_{\mu ,\beta }(s,x)=1+\sum\limits_{n=1}^{\infty
}\{x(e^{s}-1)\}^{n}\prod\limits_{k=0}^{n-1}\frac{\Gamma \lbrack \mu (k(1+%
\frac{\beta }{\mu })+\frac{\beta }{\mu })+1]}{\Gamma \lbrack \mu (k(1+\frac{%
\beta }{\mu })+\frac{\beta }{\mu }+1)+1]},  \label{eq6.3a}
\end{equation}

where $\mu $ and $\beta $ are subject to the conditions given by Eq.(\ref%
{eq1.2ab}).

If we put $x=1$ in Eq.(\ref{eq6.3}), then we immediately come to the
generating function $\mathcal{B}_{\mu ,\beta }(s)$ of the numbers $B_{\mu
,\beta }(m)$ introduced by Eq.(\ref{eq5.5})

\begin{equation}
\mathcal{B}_{\mu ,\beta }(s)=\mathcal{B}_{\mu ,\beta }(s,x)|_{x=1}=E_{\mu ,1+%
\frac{\beta }{\mu },\frac{\beta }{\mu }}(e^{s}-1).  \label{eq6.4}
\end{equation}

To get the numbers $B_{\mu ,\beta }(m)$ we differentiate $\mathcal{B}_{\mu
,\beta }(s)$ $m$ times with respect to $s$, and then set $s=0$. That is,

\begin{equation}
B_{\mu ,\beta }(m)=\frac{\partial ^{m}}{\partial s^{m}}\mathcal{B}_{\mu
,\beta }(s)|_{s=0}.  \label{eq6.5}
\end{equation}

When $\mu =1$ and $\beta =0$, using Eq.(\ref{eqA1}) we get from Eq.(\ref%
{eq5.1a})

\begin{equation*}
B_{1,0}(x,m)=e^{\text{-}x}\sum\limits_{n=0}^{\infty }n^{m}\frac{x^{n}}{n!},
\end{equation*}

which are the Bell polynomials \cite{Bell}. The new formula Eq.(\ref{eq5.5})
is a generalization of the so-called Dobi\'{n}ski relation known since 1877 
\cite{Dobinski} for the Bell numbers\footnote{%
The number of ways a set of $m$ elements can be partitioned into non-empty
subsets.}

\begin{equation*}
B_{_{1,0}}(m)=e^{\text{-}1}\sum\limits_{n=0}^{\infty }\frac{n^{m}}{n!}.
\end{equation*}

When $0<\mu \leq 1$ and $\beta =0$, using Eq.(\ref{eqA2}) we get from Eq.(%
\ref{eq5.1a})

\begin{equation*}
B_{\mu ,0}(x,m)=\sum\limits_{n=0}^{\infty }n^{m}\frac{x^{n}}{n!}E_{\mu
,1,0}^{(n)}(-x)=\sum\limits_{n=0}^{\infty }n^{m}\frac{x^{n}}{n!}E_{\mu
}^{(n)}(-x),
\end{equation*}

or in the form of a power series%
\begin{equation*}
B_{\mu ,0}(x,m)=\sum\limits_{n=0}^{\infty }n^{m}\frac{x^{n}}{n!}%
\sum\limits_{k=0}^{\infty }\frac{(n+k)!}{k!}\frac{(-x)^{k}}{\Gamma (\mu
(n+k)+1)},\quad B_{\mu ,0}(x,0)=1,
\end{equation*}

which are the fractional Bell polynomials \cite{Laskin8}. Fractional
generalization of the Dobi\'{n}ski relation is \cite{Laskin8}

\begin{equation*}
B_{\mu ,0}(m)=\sum\limits_{n=0}^{\infty }n^{m}\frac{x^{n}}{n!}E_{\mu
,1,0}^{(n)}(-1)=\sum\limits_{n=0}^{\infty }n^{m}\frac{x^{n}}{n!}E_{\mu
}^{(n)}(-1).
\end{equation*}

With $\mu +\beta =1$ and using Eq.(\ref{eqA4}) we can present Eq.(\ref%
{eq5.1a}) as

\begin{equation*}
B_{\mu ,1-\mu }(x,m)=\sum\limits_{n=0}^{\infty }n^{m}\frac{x^{n}}{n!}E_{\mu ,%
\frac{1}{\mu },\frac{1-\mu }{\mu }}^{(n)}(-x),
\end{equation*}

where $\lambda $ is defined by Eq.(\ref{eq1.171d}). The series
representation for polynomials $B_{\mu ,1-\mu }(x,m)$ has the form

\begin{equation}
B_{\mu ,1-\mu }(x,m)=\sum\limits_{n=0}^{\infty }n^{m}\frac{x^{n}}{n!}%
\sum\limits_{m=0}^{\infty }\frac{(m+n)!}{m!}(-\lambda
t)^{m}\prod\limits_{k=0}^{m+n-1}\frac{\Gamma (k+2-\mu )}{\Gamma (k+2)}.
\label{eq6.7a}
\end{equation}

\subsection{Fractional combinatorial numbers associated with polynomials $B_{%
\protect\mu ,\protect\beta }(x,m)$}

Let us introduce \textit{fractional combinatorial numbers} $S_{\mu ,\beta
}(m,l)$ associated with the polynomials $B_{\mu ,\beta }(x,m)$ using the
following equations

\begin{equation}
B_{\mu ,\beta }(x,m)=\sum\limits_{l=0}^{m}S_{\mu ,\beta }(m,l)x^{l},
\label{eq8.1}
\end{equation}

\begin{equation*}
S_{\mu ,\beta }(m,0)=\delta _{m,0},\qquad S_{\mu ,\beta }(m,l)=0,\ l=m+1,\
m+2,\ ...,
\end{equation*}

where $B_{\mu ,\beta }(x,m)$ is the polynomials given by Eq.(\ref{eq5.1a})
and the parameters $\mu $, and $\beta $ satisfy the conditions (\ref{eq1.2ab}%
).

At $x=1$, when the polynomials $B_{\mu ,\beta }(x,m)$ become the fractional
numbers, $B_{\mu ,\beta }(m)=B_{\mu ,\beta }(x,m)\left\vert _{x=1}\right. $,
Eq.(\ref{eq8.1}) goes into the expression for the fractional numbers $B_{\mu
,\beta }(m)$ in terms of fractional combinatorial numbers $S_{\mu ,\beta
}(m,l)$

\begin{equation}
B_{\mu ,\beta }(m)=\sum\limits_{l=0}^{m}S_{\mu ,\beta }(m,l).  \label{eq8.2}
\end{equation}

To find a generating function of the fractional combinatorial numbers $%
S_{\mu ,\beta }(m,l)$, let's expand the generating function $\mathcal{B}%
_{\mu ,\beta }(s,x)$ given by Eq.(\ref{eq6.1}). Upon substituting $B_{\mu
}(x,m)$ from Eq.(\ref{eq8.1}) we have the following chain of transformations

\begin{equation*}
\mathcal{B}_{\mu ,\beta }(s,x)=\sum\limits_{m=0}^{\infty }\frac{s^{m}}{m!}%
\left( \sum\limits_{l=0}^{m}S_{\mu ,\beta }(m,l)x^{l}\right) =
\end{equation*}

\begin{equation}
\sum\limits_{m=0}^{\infty }\frac{s^{m}}{m!}\left( \sum\limits_{l=0}^{\infty
}\theta (m-l)S_{\mu ,\beta }(m,l)x^{l}\right) =\sum\limits_{l=0}^{\infty
}\left( \sum\limits_{m=l}^{\infty }S_{\mu ,\beta }(m,l)\frac{s^{m}}{m!}%
\right) x^{l},  \label{eq8.3}
\end{equation}

here $\theta (l)$ is the Heaviside step function, $\theta (l)=1,\ \mathrm{%
if\ }l\geq 0$ and \ $\theta (l)=0,\ \mathrm{if\ }l<0$.

On the other hand, from Eq.(\ref{eq6.3}), we have for $\mathcal{B}_{\mu
,\beta }(s,x)$

\begin{equation}
\mathcal{B}_{\mu ,\beta }(s,x)=E_{\mu ,1+\frac{\beta }{\mu },\frac{\beta }{%
\mu }}(x(e^{s}-1))=1+  \label{eq8.4a}
\end{equation}

\begin{equation*}
\sum\limits_{n=1}^{\infty }\{x(e^{s}-1)\}^{n}\prod\limits_{k=0}^{n-1}\frac{%
\Gamma \lbrack \mu (k(1+\frac{\beta }{\mu })+\frac{\beta }{\mu })+1]}{\Gamma
\lbrack \mu (k(1+\frac{\beta }{\mu })+\frac{\beta }{\mu }+1)+1]}.
\end{equation*}

Upon comparing Eq.(\ref{eq8.3}) and Eq.(\ref{eq8.4a}), we conclude that

\begin{equation}
\sum\limits_{m=l}^{\infty }S_{\mu ,\beta }(m,l)\frac{s^{m}}{m!}%
=(e^{s}-1)^{l}\prod\limits_{k=0}^{l-1}\frac{\Gamma \lbrack \mu (k(1+\frac{%
\beta }{\mu })+\frac{\beta }{\mu })+1]}{\Gamma \lbrack \mu (k(1+\frac{\beta 
}{\mu })+\frac{\beta }{\mu }+1)+1]},\quad l=1,2,....  \label{eq8.5}
\end{equation}

Let us introduce two generating functions $\mathcal{G}_{\mu ,\beta }(s,l)$
and $\mathcal{F}_{\mu ,\beta }(s,t)$ of the fractional combinatorial numbers 
$S_{\mu ,\beta }(m,l)$,

\begin{equation}
\mathcal{G}_{\mu ,\beta }(s,l)=\sum\limits_{m=l}^{\infty }S_{\mu ,\beta
}(m,l)\frac{s^{m}}{m!}=(e^{s}-1)^{l}\prod\limits_{k=0}^{l-1}\frac{\Gamma
\lbrack \mu (k(1+\frac{\beta }{\mu })+\frac{\beta }{\mu })+1]}{\Gamma
\lbrack \mu (k(1+\frac{\beta }{\mu })+\frac{\beta }{\mu }+1)+1]},
\label{eq8.6}
\end{equation}

and

\begin{equation}
\mathcal{F}_{\mu ,\beta
}(s,t)=\sum\limits_{l=0}^{m}\sum\limits_{m=l}^{\infty }S_{\mu ,\beta }(m,l)%
\frac{s^{m}t^{l}}{m!}=E_{\mu ,1+\frac{\beta }{\mu },\frac{\beta }{\mu }%
}(t(e^{s}-1).  \label{eq8.7}
\end{equation}

Then the fractional combinatorial numbers are obtained as follows,%
\begin{equation}
S_{\mu ,\beta }(m,l)=\frac{\partial ^{m}}{\partial s^{m}}\mathcal{G}_{\mu
,\beta }(s,l)\left\vert _{s=0}\right. =  \label{eq8.8}
\end{equation}

\begin{equation*}
\prod\limits_{k=0}^{l-1}\frac{\Gamma \lbrack \mu (k(1+\frac{\beta }{\mu })+%
\frac{\beta }{\mu })+1]}{\Gamma \lbrack \mu (k(1+\frac{\beta }{\mu })+\frac{%
\beta }{\mu }+1)+1]}\frac{\partial ^{m}}{\partial s^{m}}(e^{s}-1)^{l}\left%
\vert _{s=0}\right. ,
\end{equation*}

and

\begin{equation}
S_{\mu ,\beta }(m,l)=\frac{1}{l!}\frac{\partial ^{m+l}}{\partial
s^{m}\partial t^{l}}\mathcal{F}_{\mu ,\beta }(s,t)\left\vert
_{s=0},_{t=0}\right. =  \label{eq8.9}
\end{equation}

\begin{equation*}
\frac{1}{l!}\frac{\partial ^{m+l}}{\partial s^{m}\partial t^{l}}E_{\mu ,1+%
\frac{\beta }{\mu },\frac{\beta }{\mu }}(t(e^{s}-1))|_{s=0},_{t=0},
\end{equation*}

where $l\leq m$ and $S_{\mu ,\beta }(m,0)=\delta _{m,0}$.

When $\mu =1$ and $\beta =0$ it follows from Eq.(\ref{eq8.8}) that

\begin{equation}
S_{1,0}(m,l)=\frac{1}{l!}\frac{\partial ^{m}}{\partial s^{m}}%
(e^{s}-1)^{l}\left\vert _{s=0}\right. .  \label{eq8.9a}
\end{equation}

The right side of this equation is generating function of the standard
Stirling numbers of the second kind \cite{Abramowitz}. Hence, $S_{1,0}(m,l)$
is the standard Stirling number of the second kind\footnote{%
The number of ways to partition a set of $m$ elements into $l$ non-empty
sets.} $S(m,l)\equiv S_{1,0}(m,l)$, and we conclude that

\begin{equation}
S_{\mu ,\beta }(m,l)=l!\prod\limits_{k=0}^{l-1}\frac{\Gamma \lbrack \mu (k(1+%
\frac{\beta }{\mu })+\frac{\beta }{\mu })+1]}{\Gamma \lbrack \mu (k(1+\frac{%
\beta }{\mu })+\frac{\beta }{\mu }+1)+1]}S(m,l).  \label{eq8.9b}
\end{equation}

Thus, we have represented fractional combinatorial numbers $S_{\mu ,\beta
}(m,l)$ through the the standard Stirling number of the second kind $S(m,l)$%
. It is easy to see, for example, that

\begin{equation*}
S_{\mu ,\beta }(m,1)=\frac{\Gamma (\beta +1)}{\Gamma (\mu +\beta +1)},
\end{equation*}

\begin{equation*}
S_{\mu ,\beta }(m,2)=2(2^{m-1}-1)\prod\limits_{k=0}^{1}\frac{\Gamma \lbrack
\mu (k(1+\frac{\beta }{\mu })+\frac{\beta }{\mu })+1]}{\Gamma \lbrack \mu
(k(1+\frac{\beta }{\mu })+\frac{\beta }{\mu }+1)+1]},
\end{equation*}

\begin{equation*}
S_{\mu ,\beta }(m,m-1)=\frac{m!(m-1)}{2}\prod\limits_{k=0}^{m-2}\frac{\Gamma
\lbrack \mu (k(1+\frac{\beta }{\mu })+\frac{\beta }{\mu })+1]}{\Gamma
\lbrack \mu (k(1+\frac{\beta }{\mu })+\frac{\beta }{\mu }+1)+1]},
\end{equation*}

\begin{equation*}
S_{\mu ,\beta }(m,m)=m!\prod\limits_{k=0}^{m-1}\frac{\Gamma \lbrack \mu (k(1+%
\frac{\beta }{\mu })+\frac{\beta }{\mu })+1]}{\Gamma \lbrack \mu (k(1+\frac{%
\beta }{\mu })+\frac{\beta }{\mu }+1)+1]}.
\end{equation*}

In a particular case, for $\mu =1$ and $\beta =0$, equations (\ref{eq8.6})
and (\ref{eq8.7}) turn into the well-know generating functions for standard
Stirling numbers of the second kind $S(m,l)$,

\begin{equation}
\mathcal{G}_{1,0}(s,l)=\sum\limits_{m=l}^{\infty }S(m,l)\frac{s^{m}}{m!}=%
\frac{(e^{s}-1)^{l}}{l!},\qquad l=0,1,2,...,  \label{eq8.10}
\end{equation}

and

\begin{equation}
\mathcal{F}_{1,0}(s,t)=\sum\limits_{m=0}^{\infty }\sum\limits_{l=0}^{m}S(m,l)%
\frac{s^{m}t^{l}}{m!}=\exp (t(e^{s}-1)).  \label{eq8.11}
\end{equation}

(see, Eqs.(2.17) and (2.18) in \cite{Charalambides}).

\subsubsection{Statistics of fractional probability distribution}

We now use the fractional combinatorial numbers introduced by Eq.(\ref{eq8.1}%
) to obtain the moments of the fractional non-homogeneous counting process
with the probability distribution function $P_{\mu ,\beta }(n,t)$ given by
Eqs.(\ref{eq1.20}) and (\ref{eq1.21}). By the definition of the $m^{\mathrm{%
th}}$ order moment, we have

\begin{equation}
<n_{\mu ,\beta }^{m}>=\sum\limits_{n=0}^{\infty }n^{m}P_{\mu ,\beta
}(n,t)=\sum\limits_{n=0}^{\infty }n^{m}\frac{(-z)^{n}}{n!}\frac{d^{n}}{dz^{n}%
}\left( E_{\mu ,1+\frac{\beta }{\mu },\frac{\beta }{\mu }}(z)\right)
|_{z=-\lambda _{\mu +\beta }t^{\mu +\beta }}.  \label{eq8.12}
\end{equation}

It is easy to see that $<n_{\mu ,\beta }^{m}>$ is the fractional polynomial $%
B_{\mu ,\beta }(\lambda _{\mu +\beta }t^{\mu +\beta },m)$ defined by Eq.(\ref%
{eq5.1a}),

\begin{equation*}
<n_{\mu ,\beta }^{m}>=B_{\mu ,\beta }(\lambda _{\mu +\beta }t^{\mu +\beta
},m).
\end{equation*}
Then using the definition Eq.(\ref{eq8.1}) we obtain

\begin{equation}
<n_{\mu ,\beta }^{m}>=\sum\limits_{l=0}^{m}S_{\mu ,\beta }(m,l)(\lambda
_{\mu +\beta }t^{\mu +\beta })^{l},  \label{eq8.13}
\end{equation}

where $S_{\mu ,\beta }(m,l)$ is fractional combinatorial numbers defined by
Eq.(\ref{eq8.9b}). Therefore, the fractional combinatorial numbers $S_{\mu
,\beta }(m,l)$ naturally appear in the power series over $\lambda _{\mu
+\beta }t^{\mu +\beta }$ for the $m^{\mathrm{th}}$ order moment of the
fractional non-homogeneous counting process introduced by Eqs.(\ref{eq1.1})
- (\ref{eq1.2ab}).

Using the analytical expression given by Eq.(\ref{eq8.9b}), let's list first
four moments of the fractional non-homogeneous probability distribution $%
P_{\mu ,\beta }(n,t)$. The first and second moments are given by Eqs.(\ref%
{eq1.34}) and (\ref{eq1.35}) respectively. Third order moment $<n_{\mu
,\beta }^{3}>$ is

\begin{equation}
<n_{\mu ,\beta }^{3}>=\sum\limits_{n=0}^{\infty }n^{3}P_{\mu ,\beta }(n,t)=
\label{eq823}
\end{equation}

\begin{equation*}
6\frac{\Gamma (\beta +1)}{\Gamma (\mu +\beta +1)}\frac{\Gamma (\mu +2\beta
+1)}{\Gamma (2\mu +2\beta +1)}\frac{\Gamma (2\mu +3\beta +1)}{\Gamma (3\mu
+3\beta +1)}(\lambda _{\mu +\beta }t^{\mu +\beta })^{3}+
\end{equation*}

\begin{equation*}
6\frac{\Gamma (\beta +1)}{\Gamma (\mu +\beta +1)}\frac{\Gamma (\mu +2\beta
+1)}{\Gamma (2\mu +2\beta +1)}(\lambda _{\mu +\beta }t^{\mu +\beta })^{2}+%
\frac{\Gamma (\beta +1)}{\Gamma (\mu +\beta +1)}\lambda _{\mu +\beta }t^{\mu
+\beta },
\end{equation*}

and fourth order moment $<n_{\mu ,\beta }^{4}>$ is

\begin{equation*}
<n_{\mu ,\beta }^{4}>=\sum\limits_{n=0}^{\infty }n^{4}P_{\mu ,\beta }(n,t)=
\end{equation*}

\begin{equation*}
24\frac{\Gamma (\beta +1)}{\Gamma (\mu +\beta +1)}\frac{\Gamma (\mu +2\beta
+1)}{\Gamma (2\mu +2\beta +1)}\frac{\Gamma (2\mu +3\beta +1)}{\Gamma (3\mu
+3\beta +1)}\frac{\Gamma (3\mu +4\beta +1)}{\Gamma (4\mu +4\beta +1)}%
(\lambda _{\mu +\beta }t^{\mu +\beta })^{4}+
\end{equation*}

\begin{equation}
36\frac{\Gamma (\beta +1)}{\Gamma (\mu +\beta +1)}\frac{\Gamma (\mu +2\beta
+1)}{\Gamma (2\mu +2\beta +1)}\frac{\Gamma (2\mu +3\beta +1)}{\Gamma (3\mu
+3\beta +1)}(\lambda _{\mu +\beta }t^{\mu +\beta })^{3}+  \label{eq834}
\end{equation}

\begin{equation*}
14\frac{\Gamma (\beta +1)}{\Gamma (\mu +\beta +1)}\frac{\Gamma (\mu +2\beta
+1)}{\Gamma (2\mu +2\beta +1)}(\lambda _{\mu +\beta }t^{\mu +\beta })^{2}+%
\frac{\Gamma (\beta +1)}{\Gamma (\mu +\beta +1)}\lambda _{\mu +\beta }t^{\mu
+\beta }.
\end{equation*}

When $0<\mu \leq 1$ and $\beta =0$, Eqs.(\ref{eq1.34}), (\ref{eq1.35}), (\ref%
{eq823}) and (\ref{eq834}) turn into equations for the first four moments of
the fractional Poisson process \cite{Laskin8}.

In the case when $\mu =1$ and $\beta =0$, Eqs.(\ref{eq1.34}), (\ref{eq1.35}%
), (\ref{eq823}) and (\ref{eq834}) become the well-know equations for
moments of the standard Poisson probability distribution with the parameter $%
\lambda _{1}=\lambda _{\mu +\beta }|_{\mu =1,\beta =0}$ (see, Eqs.(22) -
(24) in Ref. \cite{Mathworld}).

When $\mu =1$ and $-1<\beta \leq 0$ Eqs.(\ref{eq1.34}), (\ref{eq1.35}), (\ref%
{eq823}) and (\ref{eq834}) become the equations for the moments of the
process, defined by Eq.(\ref{eq1}). Taking into account Eq.(\ref{eq1.171a})
we come to the known equations for the moments of the standard Poisson
probability distribution with the parameter $\overline{\lambda }_{\sigma }$.

\subsection{Variance, skewness and kurtosis excess of the fractional
non-homogeneous counting process}

To find analytical expressions for the variance, skewness and kurtosis
excess of a fractional non-homogeneous counting process determined by the
probability distribution function $P_{\mu ,\beta }(n,t)$ given by Eq.(\ref%
{eq1.26}), let's introduce the central moment of the $m^{\mathrm{th}}$ order 
$\mathrm{M}_{\mu ,\beta }(m)$,

\begin{equation*}
\mathrm{M}_{\mu ,\beta }(m)=<(n-<n_{\mu ,\beta
}>)^{m}>=\sum\limits_{n=0}^{\infty }(n-<n_{\mu ,\beta }>)^{m}P_{\mu ,\beta
}(n,t)=
\end{equation*}

\begin{equation}
\sum\limits_{n=0}^{\infty }\sum\limits_{r=0}^{m}(-1)^{m-r}\binom{m}{r}%
n^{r}(<n_{\mu ,\beta }>)^{m-r}P_{\mu ,\beta }(n,t)=  \label{eq87a}
\end{equation}

\begin{equation*}
\sum\limits_{r=0}^{m}(-1)^{m-r}\binom{m}{r}(<n_{\mu ,\beta
}>)^{m-r}\sum\limits_{l=0}^{r}S_{\mu ,\beta }(r,l)(\lambda _{\mu +\beta
}t^{\mu +\beta })^{l},
\end{equation*}

where $<n_{\mu ,\beta }>$\ is the first order moment defined by Eq.(\ref%
{eq1.34}) and $S_{\mu ,\beta }(r,l)$ is given by Eq.(\ref{eq8.9b}). The
first four central moments $\mathrm{M}_{\mu ,\beta }(m)$, $m=1,2,3,4$, are
presented in the Appendix B. In terms of the central moments $\mathrm{M}%
_{\mu ,\beta }(m)$, $m=1,2,3,4$, the variance $Var_{\mu ,\beta }$, skewness $%
s_{\mu ,\beta }$, and kurtosis excess $k_{\mu ,\beta }$ of the fractional
non-homogeneous counting process are

\begin{equation}
Var_{\mu ,\beta }=\mathrm{M}_{\mu ,\beta }(2),  \label{eq87b}
\end{equation}

\begin{equation}
{\large s}_{\mu ,\beta }=\frac{\mathrm{M}_{\mu ,\beta }(3)}{\mathrm{M}_{\mu
,\beta }^{3/2}(2)},  \label{eq87c}
\end{equation}

\begin{equation}
{\large k}_{\mu ,\beta }=\frac{\mathrm{M}_{\mu ,\beta }(4)}{\mathrm{M}_{\mu
,\beta }^{2}(2)}-3.  \label{eq87d}
\end{equation}

When $0<\mu <1$ and $\beta =0$, Eqs.(\ref{eq87b}) - (\ref{eq87d}) turn into
equations for the variance, skewness and kurtosis excess of the fractional
Poisson process \cite{Laskin8}. In the case where $\mu =1$ and $\beta =0$,
Eqs.(\ref{eq87b}) - (\ref{eq87d}) become equations for the variance,
skewness and kurtosis excess of the standard Poisson probability
distribution, see, for example, Eqs.(25) - (27) in Ref. \cite{Mathworld}.
When $\mu =1$ and $-1<\beta \leq 0$ the variance, skewness and kurtosis
excess of the process defined by Eq.(\ref{eq1}) are exactly the same as for
the Poisson process if we take into account that for the process defined by
Eq.(\ref{eq1}) we have

\begin{equation*}
<n_{1,\beta }>=<n_{1+\beta }>|_{\mu =1}=\frac{\lambda _{1+\beta }t^{1+\beta }%
}{1+\beta },
\end{equation*}

which can be expressed in \textit{Poisson-like form} using Eq.(\ref{eq1.171a}%
)

\begin{equation*}
<n_{1+\beta }>\equiv <n_{\sigma }>=\overline{\lambda }_{\sigma }t^{\sigma }.
\end{equation*}

\subsubsection{New representation of the Kilbas and Saigo function}

To obtain a new representation for the Kilbas and Saigo function $E_{\mu ,1+%
\frac{\beta }{\mu },\frac{\beta }{\mu }}(z)$ defined by Eq.(\ref{eq1.16d})
we use the generating function of the Stirling numbers of the first kind%
\footnote{%
The Stirling numbers of the first kind are defined as the coefficients $%
s(m,l)$ in the expansion
\par
\begin{equation}
\frac{x!}{(x-m)!}=x(x-1)...(x-m+1)=\sum\limits_{l=0}^{m}s(m,l)x^{l}.
\label{eqML1.1}
\end{equation}%
\par
{}}

\begin{equation}
(1+z)^{n}=\sum\limits_{m=0}^{\infty }\frac{z^{m}}{m!}\sum%
\limits_{l=0}^{m}s(m,l)n^{l},  \label{eqML1}
\end{equation}

where $s(m,l)$ stands for the Stirling numbers of the first kind, see Sec.
24.1.3 in \cite{Abramowitz}. Then, we evaluate expectations of both sides of
Eq.(\ref{eqML1}) using the probability distribution function of fractional
non-homogeneous counting process $P_{\mu ,\beta }(n,\rho )$ given by Eq.(\ref%
{eq1.25})

\begin{equation}
\sum\limits_{n=0}^{\infty }P_{\mu ,\beta
}(n,1)(1+z)^{n}=\sum\limits_{n=0}^{\infty }P_{\mu ,\beta
}(n,1)\sum\limits_{m=0}^{\infty }\frac{z^{m}}{m!}\sum%
\limits_{l=0}^{m}s(m,l)n^{l}.  \label{eqML2}
\end{equation}

It is easy to see that the calculations on the left side of Eq.(\ref{eqML2})
lead to $E_{\mu ,1+\frac{\beta }{\mu },\frac{\beta }{\mu }}(z)$, and we

\begin{equation*}
E_{\mu ,1+\frac{\beta }{\mu },\frac{\beta }{\mu }}(z)=\sum\limits_{n=0}^{%
\infty }P_{\mu }(n,1)\sum\limits_{m=0}^{\infty }\frac{z^{m}}{m!}%
\sum\limits_{l=0}^{m}s(m,l)n^{l}.
\end{equation*}

Then, using the definition of the fractional numbers $B_{\mu ,\beta }(m)$
given by Eq.(\ref{eq5.5}) we obtain new representations of the Kilbas and
Saigo function

\begin{equation}
E_{\mu ,1+\frac{\beta }{\mu },\frac{\beta }{\mu }}(z)=\sum\limits_{m=0}^{%
\infty }\frac{z^{m}}{m!}\sum\limits_{l=0}^{m}s(m,l)B_{\mu ,\beta }(l)
\label{eqML2a}
\end{equation}

through Stirling numbers of the first kind and the fractional numbers $%
B_{\mu ,\beta }(l)$. Using Eq.(\ref{eq8.2}) the Kilbas and Saigo function
can be expressed through Stirling numbers of the first kind and fractional
combinatorial numbers $S_{\mu ,\beta }(m,l)$

\begin{equation}
E_{\mu ,1+\frac{\beta }{\mu },\frac{\beta }{\mu }}(z)=\sum\limits_{m=0}^{%
\infty }\frac{z^{m}}{m!}\sum\limits_{l=0}^{m}\sum\limits_{r=0}^{l}s(m,l)S_{%
\mu ,\beta }(l,r).  \label{eqML2b}
\end{equation}

The two equations (\ref{eqML2a}) and (\ref{eqML2b})\ are new representations
of the Kilbas and Saigo function.

From Eq.(\ref{eqML2a})\ we see that the Kilbas and Saigo function can be
considered as generating function of the sum $\sum\limits_{l=0}^{m}s(m,l)B_{%
\mu ,\beta }(l)$, that is

\begin{equation}
\sum\limits_{l=0}^{m}s(m,l)B_{\mu ,\beta }(l)=\frac{\partial ^{m}}{\partial
z^{m}}E_{\mu ,1+\frac{\beta }{\mu },\frac{\beta }{\mu }}(z)|_{z=0}.
\label{eqML3a}
\end{equation}

Using Eq.(\ref{eq1.16d}), we obtain the identity

\begin{equation}
\sum\limits_{l=0}^{m}s(m,l)B_{\mu ,\beta }(l)=m!\prod\limits_{k=0}^{m-1}%
\frac{\Gamma \lbrack \mu (k(1+\frac{\beta }{\mu })+\frac{\beta }{\mu })+1]}{%
\Gamma \lbrack \mu (k(1+\frac{\beta }{\mu })+\frac{\beta }{\mu }+1)+1]}.
\label{eqML3b}
\end{equation}

From Eq.(\ref{eqML2b})\ we see that the Kilbas and Saigo function can be
considered as generating function of the sum $\sum\limits_{l=0}^{m}\sum%
\limits_{r=0}^{l}s(m,l)S_{\mu ,\beta }(l,r)$, that is

\begin{equation*}
\sum\limits_{l=0}^{m}\sum\limits_{r=0}^{l}s(m,l)S_{\mu ,\beta }(l,r)=\frac{%
\partial ^{m}}{\partial z^{m}}E_{\mu ,1+\frac{\beta }{\mu },\frac{\beta }{%
\mu }}(z)|_{z=0}.
\end{equation*}

Using Eq.(\ref{eq1.16d}), we obtain the identity

\begin{equation}
\sum\limits_{l=0}^{m}\sum\limits_{r=0}^{l}s(m,l)S_{\mu ,\beta
}(l,r)=m!\prod\limits_{k=0}^{m-1}\frac{\Gamma \lbrack \mu (k(1+\frac{\beta }{%
\mu })+\frac{\beta }{\mu })+1]}{\Gamma \lbrack \mu (k(1+\frac{\beta }{\mu })+%
\frac{\beta }{\mu }+1)+1]}.  \label{eqML3c}
\end{equation}

\section{Conclusions}

A new fractional non-homogeneous counting process has been introduced and
developed using the Kilbas and Saigo three-parameter generalization of the
Mittag-Leffler function. The probability distribution function of this
process reproduces for certain set of the fractality parameters the famous
Poisson and fractional Poisson probability distributions as well as the
probability distribution function of a counting process displaying stretched
exponential interarrival times distribution.

Applications of the developed fractional non-homogeneous counting process
cover fractional compound process, the generalization of combinatorial
polynomials and numbers and their generating functions, variance, skewness
and kurtosis excess of the fractional non-homogeneous counting process
expressed through four central moments and a new representation of the
Kilbas and Saigo function.

The presented results will stimulate further development of the fundamentals
and various applications in a new field of probability theory and stochastic
processes, which can be called fractional stochastic processes and their
applications.

\section{Appendix A}

\subsection{Special cases of the Kilbas and Saigo function $E_{\protect\mu %
,1+\frac{\protect\beta }{\protect\mu },\frac{\protect\beta }{\protect\mu }%
}(z)$}

Defined by Eq.(\ref{eq1.16d}) the Kilbas and Saigo function has the
following special cases.

$\bullet $ When $\mu =1$ and $\beta =0$, Eq.(\ref{eq1.16d}) goes to

\begin{equation*}
E_{1,1,0}(z)=1+\sum\limits_{n=1}^{\infty }z^{n}\prod\limits_{k=0}^{n-1}\frac{%
\Gamma \lbrack k+1]}{\Gamma \lbrack (k+1)+1]},
\end{equation*}

Then due to

\begin{equation}
\prod\limits_{k=0}^{n-1}\frac{\Gamma \lbrack k+1]}{\Gamma \lbrack (k+1)+1]}%
=\prod\limits_{k=0}^{n-1}\frac{1}{k+1}=\frac{1}{n!},  \label{eq1.171A}
\end{equation}

where the following property of the Gamma function was used

\begin{equation}
\Gamma (z+1)=z\Gamma (z),  \label{eq1.17aB}
\end{equation}

we obtain

\begin{equation}
E_{1,1,0}(z)=1+\sum\limits_{n=1}^{\infty }\frac{z^{n}}{n!}=\exp (z).
\label{eqA1}
\end{equation}

$\bullet $ When $0<\mu \leq 1$ and $\beta =0$, Eq.(\ref{eq1.16d}) goes to

\begin{equation*}
E_{\mu ,1,0}(z)=1+\sum\limits_{n=1}^{\infty }z^{n}\prod\limits_{k=0}^{n-1}%
\frac{\Gamma \lbrack \mu k+1]}{\Gamma \lbrack \mu (k+1)+1]}.
\end{equation*}

Then due to

\begin{equation}
\prod\limits_{k=0}^{n-1}\frac{\Gamma \lbrack \mu k+1]}{\Gamma \lbrack \mu
(k+1)+1]}=\frac{1}{\Gamma (\mu n+1)},  \label{eq1.172b}
\end{equation}

we obtain

\begin{equation}
E_{\mu ,1,0}(z)=1+\sum\limits_{n=1}^{\infty }\frac{z^{n}}{\Gamma (\mu n+1)}%
=E_{\mu }(z),  \label{eqA2}
\end{equation}

here $E_{\mu }(z)$ is the Mittag-Leffler function defined by Eq.(\ref%
{eq1.172d}).

$\bullet $ When $\mu =1$ and $-1<\beta \leq 0$, Eq.(\ref{eq1.16d}) goes to

\begin{equation*}
E_{1,1+\beta ,\beta }(z)=1+\sum\limits_{n=1}^{\infty
}z^{n}\prod\limits_{k=0}^{n-1}\frac{\Gamma \lbrack \mu (k(1+\beta )+\beta
)+1]}{\Gamma \lbrack \mu (k(1+\beta )+\beta +1)+1]}.
\end{equation*}

Then due to

\begin{equation}
\prod\limits_{k=0}^{n-1}\frac{\Gamma \lbrack (k(1+\beta )+\beta )+1]}{\Gamma
\lbrack (k(1+\beta )+\beta +1)+1]}=\prod\limits_{k=0}^{n-1}\frac{1}{%
(k+1)(\beta +1)}=\frac{1}{(1+\beta )^{n}}\frac{1}{n!},  \label{eq1.17b}
\end{equation}

we obtain

\begin{equation}
E_{1,1+\beta ,\beta }(z)=1+\sum\limits_{n=1}^{\infty }\frac{1}{n!}\frac{z^{n}%
}{(1+\beta )^{n}}=\exp (\frac{z}{1+\beta }).  \label{eqA3}
\end{equation}

$\bullet $ When $\mu +\beta =1$, Eq.(\ref{eq1.16d}) goes to

\begin{equation}
E_{\mu ,\frac{1}{\mu },\frac{1-\mu }{\mu }}(z)=1+\sum\limits_{n=1}^{\infty
}z^{n}\prod\limits_{k=0}^{n-1}\frac{\Gamma (k+2-\mu )}{\Gamma (k+2)}.
\label{eqA4}
\end{equation}

\section{Appendix B}

\subsection{Four first central moments}

Using Eq.(\ref{eq87a}), the first four central moments can be calculated. We
expressed the results of the calculation in the form of a power series over
the first-order moment $<n_{\mu ,\beta }>$:

\begin{equation}
\mathrm{M}_{\mu ,\beta }(1)=0,  \label{eq88B}
\end{equation}

\begin{equation}
\mathrm{M}_{\mu ,\beta }(2)=\left( \frac{2\Gamma (\mu +\beta +1)\Gamma (\mu
+2\beta +1)}{\Gamma (\beta +1)\Gamma (2\mu +2\beta +1)}-1\right) <n_{\mu
,\beta }>^{2}+<n_{\mu ,\beta }>,  \label{eq89B}
\end{equation}

\begin{equation*}
\mathrm{M}_{\mu ,\beta }(3)=2{\LARGE (}\frac{3\Gamma (\mu +\beta +1)}{\Gamma
(\beta +1)}\frac{\Gamma (\mu +2\beta +1)}{\Gamma (2\mu +2\beta +1)}\frac{%
\Gamma (\mu +\beta +1)}{\Gamma (\beta +1)}\frac{\Gamma (2\mu +3\beta +1)}{%
\Gamma (3\mu +3\beta +1)}-
\end{equation*}

\begin{equation}
\frac{3\Gamma (\mu +\beta +1)}{\Gamma (\beta +1)}\frac{\Gamma (\mu +2\beta
+1)}{\Gamma (2\mu +2\beta +1)}+1{\LARGE )}<n_{\mu ,\beta }>^{3}+
\label{eq90B}
\end{equation}

\begin{equation*}
3\left( \frac{2\Gamma (\mu +\beta +1)}{\Gamma (\beta +1)}\frac{\Gamma (\mu
+2\beta +1)}{\Gamma (2\mu +2\beta +1)}-1\right) <n_{\mu ,\beta
}>^{2}+<n_{\mu ,\beta }>,
\end{equation*}

and

\begin{equation*}
\mathrm{M}_{\mu ,\beta }(4)=3{\LARGE (}8\frac{\Gamma ^{3}(\mu +\beta +1)}{%
\Gamma ^{3}(\beta +1)}\frac{\Gamma (\mu +2\beta +1)}{\Gamma (2\mu +2\beta +1)%
}\frac{\Gamma (2\mu +3\beta +1)}{\Gamma (3\mu +3\beta +1)}\frac{\Gamma (3\mu
+4\beta +1)}{\Gamma (4\mu +4\beta +1)}-
\end{equation*}

\begin{equation*}
8\frac{\Gamma ^{2}(\mu +\beta +1)}{\Gamma ^{2}(\beta +1)}\frac{\Gamma (\mu
+2\beta +1)}{\Gamma (2\mu +2\beta +1)}\frac{\Gamma (2\mu +3\beta +1)}{\Gamma
(3\mu +3\beta +1)}+
\end{equation*}

\begin{equation}
4\frac{\Gamma (\mu +\beta +1)}{\Gamma (\beta +1)}\frac{\Gamma (\mu +2\beta
+1)}{\Gamma (2\mu +2\beta +1)}-1{\LARGE )}<n_{\mu ,\beta }>^{4}+
\label{eq91B}
\end{equation}

\begin{equation*}
6(\frac{6\Gamma ^{2}(\mu +\beta +1)}{\Gamma ^{2}(\beta +1)}\frac{\Gamma (\mu
+2\beta +1)}{\Gamma (2\mu +2\beta +1)}\frac{\Gamma (2\mu +3\beta +1)}{\Gamma
(3\mu +3\beta +1)}-
\end{equation*}

\begin{equation*}
4\frac{\Gamma (\mu +\beta +1)}{\Gamma (\beta +1)}\frac{\Gamma (\mu +2\beta
+1)}{\Gamma (2\mu +2\beta +1)}+1)<n_{\mu ,\beta }>^{3}+
\end{equation*}

\begin{equation*}
2\left( 7\frac{\Gamma (\mu +\beta +1)}{\Gamma (\beta +1)}\frac{\Gamma (\mu
+2\beta +1)}{\Gamma (2\mu +2\beta +1)}-2\right) <n_{\mu ,\beta
}>^{2}+<n_{\mu ,\beta }>.
\end{equation*}

\end{document}